\documentclass{amsart}

\usepackage[utf8]{inputenc}
\usepackage[T1]{fontenc}

\usepackage{amsmath, amssymb, amsthm, mathtools}
\usepackage{mathrsfs, mathpazo}
\usepackage{physics}
\usepackage{geometry}

\usepackage{float}
\usepackage[all]{xy}
\usepackage{graphicx}
\usepackage{color}
\usepackage{tikz}
\usepackage{tikz-cd}
\usepackage{xparse}
\usetikzlibrary{matrix,arrows,decorations.pathmorphing}

\usepackage{cite}
\usepackage{fancyhdr}
\usepackage{hyperref}
\usepackage{cleveref}

\usepackage{indentfirst}
\usepackage{enumerate}
\usepackage{enumitem}
\usepackage{wrapfig}
\usepackage{adjustbox}
\usepackage{leftidx}
\usepackage{stackrel}
\usepackage{comment}

\newtheorem*{theorem*}{Theorem}
\newtheorem{example}{Example}[section]
\newtheorem{theorem}{Theorem}[section]
\newtheorem{proposition}{Proposition}[section]
\newtheorem{corollary}{Corollary}[section]
\newtheorem{definition}{Definition}[section]
\newtheorem{lemma}{Lemma}[section]
\newtheorem{remark}{Remark}[section]
\newtheorem{question}{Question}[section]

\crefname{theorem}{Theorem}{theorem}
\crefname{lemma}{Lemma}{lemma}
\crefname{remark}{Remark}{remark}
\crefname{corollary}{Corollary}{corollary}
\crefname{proposition}{Proposition}{proposition}
\crefname{example}{Example}{example}
\crefname{definition}{Definition}{definition}
\crefname{notation}{Notation}{notation}
\crefname{appendix}{Appendix}{appendix}
\crefname{section}{Section}{section}
\crefname{question}{Question}{question}

\newtheorem{innercustomthm}{Theorem}
\newenvironment{nnumthm}[1]
  {\renewcommand\theinnercustomthm{#1}\innercustomthm}
  {\endinnercustomthm}

\newcommand{\R}{\mathbb R}  
\newcommand{\N}{\mathbb N} 

\newcommand{\FA}{\mathfrak A}
\newcommand{\FM}{\mathfrak M}
\newcommand{\FN}{\mathfrak N}
\newcommand{\FR}{\mathfrak R}
\newcommand{\FH}{\mathfrak H}

\newcommand{\CB}{\mathcal B}

\newcommand{\CH}{\mathcal H} 

\newcommand{\CK}{\mathcal K}

\newcommand{\CR}{\mathcal R}

\DeclareMathOperator{\Ran}{ran}
\DeclareMathOperator{\Ind}{Ind}

\begin{document}
\title{On cleanness of von Neumann algebras}

\author{Lu Cui}
\address{Institute of Mathematics, Academy of Mathematics and Systems Science, Chinese Academy of Sciences, Beijing 100190, China}
\email{cuilu17@mails.ucas.ac.cn}
\thanks{}

\author{LinZhe Huang}
\address{Yau Mathematical Sciences Center, Tsinghua University, Beijing 100084, China}
\email{huanglinzhe@tsinghua.edu.cn}
\thanks{}

\author{Wenming Wu}
\address{School of Mathematical Sciences, Chongqing Normal University, Chongqing 401331, China}
\email{wuwm@amss.ac.cn}
\thanks{Wu was supported in part by NFS of China (No. 11871127, No. 11971463) and Chongqing Science and Technology Commission (No. cstc2019jcyj-msxmX0256).}

\author{Wei Yuan}
\address{Institute of Mathematics, Academy of Mathematics and Systems Science, Chinese Academy of Sciences, Beijing 100190, China}
\address{School of Mathematical Sciences, University of Chinese Academy of Sciences, Beijing 100049, China}
\email{wyuan@math.ac.cn}
\thanks{Yuan was supported in part by NFS of China (No. 11971463, No. 11871303, No. 11871127).}

\author{Hanbin Zhang}
\address{School of Mathematics (Zhuhai), Sun Yat-sen University, Zhuhai 519082, Guangdong, China}
\email{zhanghb68@mail.sysu.edu.cn}
\thanks{Zhang was supported by NSF of China (No. 11901563).}

\subjclass[2010]{Primary 47C15, Secondary 16U99, 46L10, 47A65}
\keywords{Von Neumann algebras; Clean rings; Idempotents; Projections}

\begin{abstract}
    A unital ring is called \textit{clean} (resp. \textit{strongly clean}) if every element can be written as the sum of an invertible element and an idempotent (resp. an invertible element and an idempotent that commutes). T.Y. Lam proposed a question: which von Neumann algebras are clean as rings? In this paper, we characterize strongly clean von Neumann algebras and prove that all finite von Neumann algebras and all separable infinite factors are clean.
\end{abstract}

\maketitle

\section{Introduction}
In his study of continuous geometry \cite{vN60}, in order to generalize the classical coordinatization theorem, von Neumann introduced a new class of rings, which he called \textit{regular rings}; see \cite{vN36}. Let $\FR$ be a unital ring. An element $T\in \FR$ is called (von Neumann) \textit{regular} if there exists an \textit{inner inverse} $U \in \FR$ for $T$, i.e., $TUT = T$. Then $\FR$ is called (von Neumann) \textit{regular} if every element of $\FR$ is regular. Moreover, if $T$ has an invertible inner inverse, then we call $T$ \textit{unit-regular} and $\FR$ is called \textit{unit-regular} if every element of $\FR$ is unit-regular; see \cite{GE68}. Studies related to regular rings have been one of the central topics in ring theory; see \cite{G91} for a monograph on regular rings.

An important generalization of regular rings arose in the study of direct sum decomposition theory. Crawley and J\'onsson \cite{CJ64} introduced the exchange property for modules to study the isomorphic refinement problem. A module $\FA$ is said to have the exchange property if for any module $\FH$ and any two decompositions $\FH = \FA \oplus \FM = \sum_{i\in I} \FN_i$, there exist submodules $\FN_i' \subseteq \FN_i$ such that $\FH = \FA \oplus \left (\sum_{i\in I} \FN_i'\right)$. Based on the work of Crawley and J\'onsson \cite{CJ64}, Warfield introduced the notion of \textit{exchange rings} in \cite{WJ72}. An exchange ring is a ring $\FR$ such that the left module $_{\FR}\FR$ (or equivalently \cite[Corollary 2]{WJ72}, the right module $\FR_{\FR}$) has the exchange property. Exchange rings are not only very useful in decomposition theory, but also form a large class of rings. Besides regular rings, the class of exchange rings includes $\pi$-regular rings \cite{S86}, semiperfect rings \cite{WJ72} and (unital) C$^*$-algebras of real rank zero \cite[Theorem 7.2]{AGMP98}. We refer to \cite{KLN16,KLN17} for some recent progress on exchange rings and references therein.

Later, Nicholson \cite{N77} showed that exchange rings are exactly those for which idempotents can be lifted modulo every left ideal, i.e., for every ring element $A$ with $A- A^2$ in a left ideal $\mathcal{L}$, there exists an idempotent $P$ in the ring such that $P - A \in \mathcal{L}$. In addition, in order to characterize those exchange rings with central idempotents, Nicholson introduced a new class of rings.
\begin{definition}
A ring is called \textit{clean} if every element can be written as the sum of an invertible element and an idempotent.
\end{definition}
Clean rings can be regarded as an additive analogy of unit-regular rings, as every element in a unit-regular ring can be written as the product of an invertible element and an idempotent. It is very interesting that unit-regularity implies cleanness of rings; see \cite{CY94,CK01}. Although, in general, clean rings are not necessarily exchange rings (see \cite[Example 1]{H77}, \cite[p. 4746]{CY94}, and \cite[Example 3.1]{S12}), many natural examples of exchange rings are clean rings, including abelian exchange rings \cite{N77}, strongly $\pi$-regular rings \cite{N99} and endomorphism rings of continuous modules \cite{CKLNZ06}. Among them, abelian exchange rings and strongly $\pi$-regular rings satisfy the following stronger condition: every element can be written as the sum of an invertible element and an idempotent that commutes. Such rings are called \textit{strongly clean rings}. We refer to \cite{NZ05} for a survey on clean rings. The recent paper \cite{NS18} provided a very nice discussion on connections between (unit-)regularity and (strong) cleanness of rings. Also see \cite{KLNS20} for a very recent element-wise study on clean rings. Moreover, Aghajani and Tarizadeh \cite{AT20} provided some new characterizations of clean rings from a topological viewpoint; Bezhanishvili, Morandi and Olberding \cite{BMO13}, based on Gelfand-Neumark-Stone duality, found that the notion of clean ring appears naturally in a study of category theory (the category of bounded Archimedean lattice-ordered algebras over $\mathbb R$).

In \cite{AGMP98}, Ara et al. proved that a unital C$^*$-algebra is an exchange ring if and only if it has real rank zero. Based on this result, together with a result of Lin \cite{L96}, Ara et al. \cite{AGMR20} partially verified a conjecture of Zhang on the topological $K_1$-group of any unital C$^*$-algebra with real rank zero (see \cite{Z91} and references therein). Moreover, it is known \cite[Proposition 1.3]{BP91} that every von Neumann algebra has real rank zero. Consequently, every von Neumann algebra is an exchange ring. It is interesting to know which C$^*$-algebras with real rank zero are clean rings. In 2005, at the Conference on Algebra and Its Applications held at Ohio University, Athens, OH, T.Y. Lam, as recorded by \cite[Introduction]{LV10}, proposed the following more specific question. 
\begin{question}\label{ques: Lam clean}
Which von Neumann algebras are clean as rings?
\end{question}
It follows immediately from the theory of Jordan canonical forms that all finite-dimensional von Neumann algebras are clean, since every finite-dimensional von Neumann algebra is $*$-isomorphic to a finite direct sum of full matrix algebra over $\mathbb C$.
Trying to answer Question \ref{ques: Lam clean}, Va\v{s} \cite{LV10} pointed out that it is more natural to utilize the fact that a von Neumann algebra is a $*$-ring (i.e. a ring with involution) and that the projections are $*$-invariant idempotents. Therefore, she introduced the following more specified cleanness for $*$-rings.
\begin{definition}
A $*$-ring $\FR$ is called a \textit{$*$-clean ring} if every element of $\FR$ can be written as the sum of an invertible element and a projection.
\end{definition}
Recall that a ring is called \textit{almost clean} if its every element can be written as the sum of a non-zero-divisor element (neither a left nor a right zero-divisor) and an idempotent (see \cite{Mc03}). Motivated by this definition, Va\v{s} also introduced the notion of almost $*$-cleanness: a $*$-ring is \textit{almost $*$-clean} if its every element can be written as the sum of a non-zero-divisor and a projection. Towards Question \ref{ques: Lam clean}, Va\v{s} \cite[Corollary 14]{LV10} proved that all finite type I von Neumann algebras are almost $*$-clean.  Later, Akalan and Va\v{s} \cite[Corollary 3.10]{AV13} showed that all finite von Neumann algebras are almost clean. Similar to strongly clean ring, a \textit{strongly $*$-clean ring} is a $*$-ring such that every element can be written as the sum of an invertible element and a projection that commutes. The study of $*$-cleanness of rings has received a lot of attention and has recently found applications in coding theory; see \cite{HZ21} and references therein.

It is known that every finite-dimensional von Neumann algebra is $*$-clean \cite[Proposition 4]{LV10}. While, in general, an infinite-dimensional von Neumann algebra may not be $*$-clean. In Section 2, we shall show that properly infinite von Neumann algebras are neither $*$-clean nor strongly clean. Based on this observation, we provide the following characterizations of strongly clean and strongly $*$-clean von Neumann algebras.

\begin{theorem}\label{thm:char_strongly_clean}
    Let $\mathfrak A$ be a von Neumann algebra.
    \begin{enumerate}
        \item $\mathfrak A$ is strongly clean if and only if there exists a finite number of mutually orthogonal central projections $P_i$ with sum $I$ such that $P_i \mathfrak A$ is of type $I_{n_i}$, $n_i < \infty$.

        \item $\mathfrak A$ is strongly $*$-clean if and only if $\mathfrak A$ is abelian.
    \end{enumerate}
\end{theorem}

Meanwhile, we continue the study of \cref{ques: Lam clean} and a large class of von Neumann algebras is shown to be clean. Firstly, we prove that all finite von Neumann algebras are clean. More precisely, we prove the following theorem.

\begin{theorem}\label{thm:finite_is_clean}
    Let $T$ be an operator in a finite von Neumann algebra $\mathfrak A$. Then there exists an idempotent $P \in \mathfrak A$ such that $T-P$ is invertible and $\|(T-P)^{-1}\| \leq 4$. In particular, $\mathfrak A$ is clean.
\end{theorem}

We also prove that a von Neumann algebra is almost $*$-clean if and only if it is finite (see \cref{rem:char_almost}) and that every properly infinite von Neumann algebra is not *-clean. Furthermore we point out that all operators with closed ranges in finite von Neumann algebras are $*$-clean (see \cref{rem:finite von Neumann algebra of type $I_n$ $*$-clean}). Next, we consider infinite von Neumann algebras and obtain the following theorem.

\begin{theorem}\label{thm:separable_factor_is_clean}
   Every separable infinite factor is clean.
\end{theorem}

We hope that these results can shed more light on Question \ref{ques: Lam clean}. The following sections are organized as follows. In Section \ref{Sec:strongly_clean}, we consider strong cleanness of von Neumann algebras and prove Theorem \ref{thm:char_strongly_clean}. In Section \ref{sec:aux_results}, we provide some auxiliary results to be used in later proofs. Next, in Section \ref{sec:finite_is_clean}, we consider finite von Neumann algebras and prove Theorem \ref{thm:finite_is_clean}. In Section \ref{sec:separable_factor_is_clean}, we deal with infinite factors (type I$_\infty$, II$_\infty$ or III factor) acting on a separable Hilbert space and prove \cref{thm:separable_factor_is_clean}. We end this paper with some remarks and questions in Section \ref{sec:concluding remarks}. We refer the reader to \cite{KR97,T02} for the general theory of von Neumann algebras.

\section{Strongly clean von Neumann algebras}\label{Sec:strongly_clean}

In this section, we characterize strongly clean von Neumann algebras. The proof of \cref{thm:char_strongly_clean} is divided into several lemmas. The first lemma is an observation which motivates us to prove \cref{thm:char_strongly_clean}.

\begin{lemma}\label{lem:infinite_not_strong_clean}
Let $\FA$ be a properly infinite von Neumann algebra acting on a Hilbert space $\CH$. Then $\FA$ is neither $*$-clean nor strongly clean.
\end{lemma}

\begin{proof}
    By Lemma 6.3.3 in \cite{KR97}, we can construct a countably infinite orthogonal family $\{E_n\}_{n=1}^\infty$ of projections with sum $I$ such that $E_n \sim I$. Let $S = \sum_{n=1}^{\infty} V_n$, where $V_n$ is a partial isometry such that $V^*_nV_n = E_n$ and $V_nV_n^* = E_{n+1}$. It is clear that $S$ is an isometry such that $SS^* = I-E_1$. In particular, $S$ is a semi-Fredholm operator with the Fredholm index $\Ind(S) := \dim \ker(S)- \dim \Ran(S)^{\perp} \neq 0$. By Theorem 5.22 in Chapter 4 in \cite{K95}, $3S -P$ is a semi-Fredholm operator with $\Ind(3S-P) = \Ind(S) \neq 0$ for every projection $P \in \FA$ (see also \cite{MB68, MB69}). Therefore $3S$ can not be written as the sum of an invertible operator and a projection. Consequently $\FA$ is not $*$-clean.

    Let $Q$ be an idempotent in $\FA$ such that $QS = SQ$. Note that $E_{1}S = 0$ and $E_{i+1} S = SE_{i}$ for every $i \geq 1$. We have $E_{1}QE_{j+1}= E_{1}QSE_jS^* = E_{1}SQE_jS^* = 0$ and 
    \begin{align*}
         E_{i+1}Q E_{j+1} = E_{i+1}Q(I-E_1) E_{j+1} = E_{i+1}S Q S^*E_{j+1} = S E_{i} Q E_{j} S^* \quad \forall i, j \geq 1.
    \end{align*}
In particular, $E_iQE_j = 0$ if $j > i$. This implies that $(E_1QE_1)^2 = E_1 Q E_1$ and $(E_1 - E_1 QE_1)(S-Q) = 0$. Thus $S - Q$ is not invertible if $E_1 QE_1 \neq E_1$. If $E_1 Q E_1 = E_1$, then $E_i Q E_i = SE_{i-1}QE_{i-1} S^* = \cdots = S^{i-1} E_1 Q E_1 (S^{i-1})^* = E_{i}$ for every $i > 1$. Recall that $Q^2 = Q$ and $E_iQE_j = 0$ for every $j > i$. We have $Q = I$. Since $S$ is a unilateral shift, $S-I$ is not invertible. Therefore $\FA$ is not strongly clean.
\end{proof}

\begin{lemma}\label{lem:shift_tensor_not_strong_clean}
    Let $\CH := \oplus_{n=1}^{\infty} \mathbb{C}^n$. We use $\{e_{n, i}\}_{i=1}^n$ to denote the canonical orthogonal basis of $\mathbb{C}^n$. Let $V := \oplus_{n=1}^{\infty}  V_n$, where $V_n$ is the upper shift matrix in $M_n(\mathbb{C})$, i.e., $V_n e_{n, i} = e_{n, i-1}$ ($e_{n,0} = 0$). Then for every Hilbert space $\CK$, the operator $I \otimes V$ is not strongly clean in $\CB(\CK \otimes \CH)$.
\end{lemma}

\begin{proof}
Let $W: l^2(\N) \otimes l^2(\N) \to \CH$ be a unitary operator defined as follows
\begin{align*}
    W (e_{k} \otimes e_{i}) = e_{i+k-1, i},
\end{align*}
where $\{e_k\}_{k=1}^{\infty}$ is the canonical orthogonal basis of $l^2(\N)$ and $k,i \geq 1$. Note that 
\begin{align*}
    W^* V W (e_{k} \otimes e_{i})  = W^*Ve_{i+k-1, i} = W^*e_{i+k-1, i-1}= \begin{cases}
        0 & i = 1,\\
        e_{k+1} \otimes e_{i-1} & i > 1.
    \end{cases} 
\end{align*} 
We have $W^*VW  = S \otimes S^*$, where $S: l^2(\N) \to l^2(\N)$ the unilateral shift, i.e., $S e_i = e_{i+1}$. Therefore, we only need to show that there exists no idempotent $P \in \CB(\CK \otimes l^2(\N) \otimes l^2(\N))$ such that $P(I \otimes S \otimes S^*) = (I \otimes S \otimes S^*)P$ and $I \otimes S \otimes S^* -P$ is invertible.

    Let $P \in \CB(\CK \otimes l^2(\N) \otimes l^2(\N))$ be an idempotent such that $P (I \otimes S \otimes S^*) = (I \otimes S \otimes S^*)P$. We use $\{E_{ij}\}$ to denote the canonical matrix units of $\CB(l^2(\N))$. An argument similar to the one used in the proof of \cref{lem:infinite_not_strong_clean} shows that $(I \otimes I \otimes E_{jj})P(I \otimes I \otimes E_{ii}) = 0$ for every $j > i$ and 
\begin{align}\label{equ:shift_tensor_not_strong_clean_1}
      \left(I \otimes I \otimes  E_{nn} \right)P \left(I \otimes I \otimes E_{nn} \right) = \left (I \otimes S^{n-1} \otimes E_{1n} \right)^*P \left (I \otimes S^{n-1} \otimes E_{1n} \right). 
   \end{align} 
   In particular, $(I \otimes I \otimes E_{11})P(I \otimes I \otimes E_{11}) = P(I \otimes I \otimes E_{11})$. Since $P^2 = P$, we have
\begin{align*}
    \left(I \otimes S \otimes S^* - P \right) \left (I \otimes I \otimes E_{11} \right) \left (I \otimes I \otimes I -P \right) \left (I \otimes I \otimes E_{11} \right) = 0. 
\end{align*}
    Therefore, $I \otimes S \otimes S^* - P$ is not invertible if $(I \otimes I \otimes E_{11})P(I \otimes I \otimes E_{11}) \neq I \otimes I \otimes E_{11}$. Assume that $(I \otimes I \otimes E_{11})P(I \otimes I \otimes E_{11}) = I \otimes I \otimes E_{11}$. By \cref{equ:shift_tensor_not_strong_clean_1}, $(I \otimes I \otimes  E_{nn})P(I \otimes I \otimes E_{nn}) = I \otimes I \otimes E_{nn}$ for every $n$. This implies that $P = I \otimes I \otimes I$. Note that
\begin{align*}
   \left \|\frac{1}{\sqrt{n}} \left(V_n-I_n \right) \left (\sum_{i=1}^n e_{n,i} \right) \right \| = \frac{1}{\sqrt{n}},
\end{align*}
we have $\left \| \left (V_n-I_n \right)^{-1} \right \| \geq \sqrt{n}$. Therefore $V - I$ is not invertible. This implies that $I \otimes S \otimes S^* - I \otimes I \otimes I$ is not invertible.
\end{proof}

\begin{lemma}\label{lem:type_II_1_not_strong_clean}
Let $\FA$ be a type II$_1$ von Neumann subalgebra of $\CB(\CH)$. Then $\FA$ is not strongly clean.
\end{lemma}

\begin{proof}
    Let $\CH_1:= \oplus_{n=1}^{\infty} \mathbb{C}^n$ and $V:=\oplus_{n=1}^{\infty} V_n$, where $V_n$ is the upper shift matrix in  $M_n(\mathbb{C})$, i.e., $V_n e_{n, i} = e_{n, i-1}$ ($e_{n,0} = 0$), where $\{e_{n, i}\}_{i=1}^n$ is the canonical orthogonal basis of $\mathbb{C}^n$.
 
    By Lemma 6.5.6 in \cite{KR97}, there exists a projection $E_{1, 1}$ such that $\Delta(E_{1,1}) = \Delta(I-E_{1,1}) = \frac{1}{2} I$, where $\Delta$ is the center-valued dimension function of $\FA$ (see, for example, Theorem 8.4.3 in \cite{KR97}). Invoke Lemma 6.5.6 in \cite{KR97} again, we can find two mutually equivalent orthogonal subprojections $E_{2,1}$, $E_{2,2}$ of $I- E_{1,1}$ such that 
    \begin{align*}
        \Delta(E_{2,1}) = \Delta(E_{2,2}) = \Delta(I-E_{1,1} - E_{2,1} - E_{2,2}) = \frac{1}{3!}I.
    \end{align*} 
Proceed as above, we can show that there exists a family of mutually equivalent orthogonal projections $\{E_{n,1}, \ldots, E_{n, n}\}$ in $\FA$ such that $E_{n,i} < I - \sum_{k=1}^{n-1} \sum_{j=1}^k E_{k, j}$, and
   \begin{align*}
       \Delta(E_{n,1}) = \cdots = \Delta(E_{n, n}) = \Delta \left (I - \sum_{k=1}^{n} \sum_{j=1}^k E_{k, j} \right) =\frac{1}{(n+1)!}I
   \end{align*} 
    for each $n \geq 1$. In particular, $I = \sum_{n=1}^{\infty} \sum_{i=1}^{n} E_{n,i}$ since $1 = \sum_{k=1}^{\infty} \frac{k}{(k+1)!}$.

For every $n\ge 2$, let $U_n:= \sum_{i=2}^{n} U_{n,i}$, where $U_{n,i}$ is a partial isometry in $\FA$ such that $U_{n,i}^*U_{n,i}= E_{n,i}$ and $U_{n,i}U_{n,i}^* = E_{n,i-1}$. Let
$U:= \sum_{n=2}^{\infty}U_n$. Note that for any $i\le n$ and $j\le k$ we have $\dim E_{n,i}\CH = \dim E_{k,j}\CH$, since $\CH$ is an infinite dimensional Hilbert space. Let $\CH_0$ be a Hilbert space such that $\dim \CH_0 = \dim E_{1,1}$ and $W_n: \CH_0 \to E_{n,n}\CH$ be a unitary operator for every $n \geq 1$. If $W: \CH_0 \otimes \CH_1 = \oplus_{n=1}^{\infty} \CH_0 \otimes \mathbb{C}^n \to \CH$ is the unitary operator defined as follows
    \begin{align*}
        W \left (\xi \otimes e_{n,i} \right) := U_n^{n-i} W_n\xi, \quad \xi \in \CH_0,
    \end{align*}
then $W^* U W = I_0 \otimes V$, where $I_0$ is the unit of $\CB(\CH_0)$. By \cref{lem:shift_tensor_not_strong_clean}, $I_0 \otimes V$ is not strongly clean.
\end{proof}

The following simple fact together with Schur's unitary triangularization theorem will be used in the proofs of \cref{lem:sep_spectrum_case} and \cref{lem:I_strongly_clean} to estimate the norms of inverses of matrices.
 
\begin{lemma}\label{lem:nil_invert_norm}
    Let $A, B \in M_n(\mathbb{C})$. If $A$ is invertible and $\left (A^{-1}B \right)^n = 0$, then $A - B$ is invertible and 
    \begin{align*}
        \left \| (A - B)^{-1} \right \| \leq \left \|A^{-1} \right \| \left( \sum_{k=0}^{n-1} \left \|A^{-1}B \right\|^k \right ). 
    \end{align*}
\end{lemma}

\begin{proof}
    Note that $(A-B)^{-1} = \left [\sum_{k=0}^{n-1}\left (A^{-1}B\right)^k \right]A^{-1}$. Then $\left \| (A - B)^{-1} \right \| \leq \left \|A^{-1} \right \| \left ( \sum_{k=0}^{n-1} \left \|A^{-1}B \right\|^k \right)$. 
\end{proof}

\begin{lemma}\label{lem:sep_spectrum_case}
    Given $A= (a_{ij}) \in M_n(\mathbb{C})$, $n \geq 2$. Let $C_1 = 3n \sum_{k=0}^{n-1}(8n \left \|A \right \|)^k$ and $C_2 = 4 \sum_{k=0}^{n-1}[8 \left \|A \right \| + 4C_1]^k$. If there exist $a, b \in \sigma(A)$ such that $|a| < \frac{1}{4}$ and $|b| > \frac{3}{4}$, then there exists $r \in [\frac{1}{4}, \frac{3}{4}]$ satisfying
    \begin{enumerate}
        \item $\{z: |z|=r\} \cap \sigma(A) = \emptyset$,
        \item $\left \|(zI-A)^{-1} \right \| \leq 4n \sum_{k=0}^{n-1}(8n\left \|A \right\|)^k$ for every $|z| = r$.
    \end{enumerate}
Therefore $P := \frac{1}{2\pi i}\int_{\Gamma} (zI-A)^{-1} dz$ is an idempotent such that $AP = PA$, $A-P$ is invertible, $\|P\| \leq C_1$ and $\left \|(A-P)^{-1} \right \| \leq C_2$, where $\Gamma$ is the circle, oriented counterclockwise, forming the boundary of the disk $\{z: |z| \leq r\}$.
\end{lemma}

\begin{proof}
    By Schur's unitary triangularization theorem, we may assume that $a_{ij} =0$ for every $i > j$ and $|a_{11}| \geq |a_{22}| \geq \cdots \geq |a_{nn}|$. Since $|a_{11}| > \frac{3}{4}$ and $|a_{nn}| < \frac{1}{4}$, the set $\left \{a_{ll}: \frac{1}{4} < |a_{ll}| < \frac{3}{4} \right\}$ contains at most $n-2$ elements. Then there exists $m$ such that one of the following conditions holds: 
    \begin{itemize}
        \item $|a_{m+1m+1}| \leq \frac{1}{4}$ and $|a_{mm}|-\frac{1}{4} \geq \frac{1}{2n}$,
        \item $a_{m-1m-1} \geq \frac{3}{4}$ and $\frac{3}{4} - |a_{mm}| \geq \frac{1}{2n}$,
        \item $|a_{mm}|, |a_{m+1m+1}| \in [\frac{1}{4}, \frac{3}{4}]$ and $|a_{mm}| - |a_{m+1m+1}| \geq \frac{1}{2n}$.
    \end{itemize}
    Therefore, there exists $r \in \left [\frac{1}{4}, \frac{3}{4} \right]$ such that $|r-|a_{ii}|| \geq \frac{1}{4n}$ for every $i$. 
    
Let $A_0$ be the diagonal matrix $(a_{ij} \delta_{ij}) \in M_n(\mathbb{C})$, and $A_1 := A - A_0$. Since $\|A_1\| \leq 2\|A\|$ and $|r - a_{ii} | \geq \frac{1}{4n}$, we have
    \begin{align*}
        \left \|(zI-A)^{-1} \right \| \leq \left \|(zI - A_0)^{-1} \right \| \left( \sum_{k=0}^{n-1} \left \|(zI - A_0)^{-1}A_1 \right\|^k \right ) \leq  4n \sum_{k=0}^{n-1}(8n \left \|A \right \|)^k
    \end{align*}
for every $z \in \Gamma$, i.e., $|z| = r$, by \cref{lem:nil_invert_norm}. Assume that $|a_{j+1 j+1}| < r < |a_{jj}|$. Recall that $|r| \leq \frac{3}{4}$. We have $P = (p_{ij}) := \frac{1}{2\pi i}\int_{\Gamma} (zI-A)^{-1} dz$ is an idempotent such that $\|P\| \leq C_1$, $PA = AP$, $p_{ij} = 0$ for every $i > j$, $p_{11} = \cdots = p_{jj} =0$, $p_{j+1 j+1} = \cdots = p_{nn} = 1$ (see Proposition 4.11 in Chapter VII in \cite{JC90}). Let $P_0:= (p_{ij} \delta_{ij})$ and $P_1: = P-P_0$. Note that $\|(A_0 -P_0)^{-1}\| \leq 4$ and $\|P_1 -A_1\| \leq \|P\|+ 2\|A\|$. By \cref{lem:nil_invert_norm}, we have 
\begin{align*}
    \left \|(A-P)^{-1} \right \| = \left \|(A_0 - P_0)^{-1} \right \|\sum_{k=0}^{n-1} \left\| (A_0 - P_0)^{-1} (P_1-A_1) \right \|^k \leq 4 \sum_{k=0}^{n-1}\left[8\|A\| + 4C_1 \right]^k. 
\end{align*}
\end{proof}

\begin{lemma}\label{lem:I_strongly_clean}
   Let $\FA$ be a type I$_n$ von Neumann algebra, where $n < \infty$. Then $\FA$ is strongly clean.
\end{lemma}

\begin{proof}
    We may assume $\FA = M_n(L^{\infty}(X, \mu))$, where $X$ is a locally compact space with the positive Radon measure (see Theorem 6.6.5 in \cite{KR97} and Theorem 1.18 in Chapter III in \cite{T02}). Let $A = \left (a_{ij} \right)_{n\times n} \in \FA$, where $a_{ij} \in L^{\infty}(X)$. In the following, we use $A(x)$ to denote $\left (a_{ij}(x) \right)_{n\times n}$ for a fixed element $x \in X$, $I_n$ to denote the unit matrix in $M_n(\mathbb{C})$. We may assume that $\|A(x)\| \leq \|A\|$ for every $x \in X$. Let $\{X_1, \ldots, X_{m}\}$ be a partition of $X$ such that $\left \|A(x)-A(y) \right \| \leq \frac{1}{2C_2}$ for every $x,y \in X_k$, $k = 1, \ldots, m$, where $C_2 = 4 \sum_{k=0}^{n-1}\left [8\|A\| + 4C_1 \right ]^k$ and $C_1 = 3n \sum_{k=0}^{n-1}(8n\|A\|)^k$. Note that $M_n(L^\infty(X)) \cong \oplus_{k=1}^{m} M_n(L^\infty(X_k))$. By considering $(a_{ij}(x)|_{X_k})_{n\times n} \in M_n(L^\infty(X_k))$, we may assume that $\|A(x)-A(y)\| \leq \frac{1}{2C_2}$ for every $x, y \in X$.

    Assume that there exists $x \in X$ such that $\sigma(A(x)) \subset \left \{z: |z|\leq \frac{3}{4} \right \}$. By Schur's unitary triangularization theorem and \cref{lem:nil_invert_norm}, we have 
    \begin{align*}
        \left \|(A(x)-I_n)^{-1} \right \| \leq  4\sum_{k=0}^{n-1}(8\|A\|)^k \leq C_2. 
    \end{align*} 
    Since $\left \|(A(x) - I_n)^{-1}(A(x) - A(y) \right \| \leq \frac{1}{2}$, $A(y) - I_n = (A(x) - I_n)[I_n - (A(x) - I_n)^{-1}(A(x) - A(y))]$ is invertible and $\|(A(y) - I_n)^{-1}\| \leq 2C_2$. Therefore $A - I$ is invertible. Similarly, if there exists $x \in X$ such that $\sigma(A(x)) \subset \{z: |z| \geq \frac{1}{4}\}$, then $A(y) = A(x)\left [I_n - A(x)^{-1}(A(x) - A(y)) \right]$ is invertible and $\left \|A(y)^{-1} \right\| \leq 2C_2$. Consequently, $A$ is invertible.

    Assume that $\sigma(A(y)) \cap \left \{z: |z| < \frac{1}{4} \right \} \neq \emptyset$ and $\sigma(A(y)) \cap \left \{z: |z| > \frac{3}{4} \right \} \neq \emptyset$ for every $y \in X$. Let $x \in X$. By \cref{lem:sep_spectrum_case}, there exists $r \in \left [\frac{1}{4}, \frac{3}{4} \right ]$ such that $\left \|(zI_n - A(x))^{-1} \right \| \leq C_2$ for every $|z|= r$. Therefore $zI_n - A(y) = (zI_n- A(x)) \left [I_n - (zI_n-A(x))^{-1}(A(y) - A(x)) \right]$ is invertible and $\left \|(z I_n - A(y))^{-1} \right \| \leq 2C_2$ for every $y \in X$ and $|z|= r$. Recall that $(zI_n - A(y))^{-1}$ equals $\frac{1}{\det \left (zI_n - A(y) \right )}$ times the cofactor matrix of $zI_n - A(y)$. We have $(zI-A)^{-1}(y) = (zI_n - A(y))^{-1}$ is in $M_n(L^\infty(X))$ for every $|z| = r$. In particular, $zI - A$ is invertible and $\left \|(zI - A)^{-1} \right \| \leq 2C_2$. By Proposition 4.11 in Chapter VII in \cite{JC90}, $P(y) := \frac{1}{2\pi i}\int_{|z|=r} (zI_n- A(y))^{-1} dz$ is an idempotent in $M_n(L^\infty(X))$ such that $PA = AP$. Since $r \leq \frac{3}{4}$, $\|P(y)\| \leq \frac{3C_2}{2}$. By Schur's unitary triangularization theorem and \cref{lem:nil_invert_norm}, we have $\left \|(A(y) - P(y))^{-1} \right \| \leq 4 \sum_{k=0}^{n-1}\left [8\|A\| + 6C_2 \right]^k$. Thus $A-P$ is invertible. 

\end{proof}

We are now ready to prove \cref{thm:char_strongly_clean}. The second part of \cref{thm:char_strongly_clean} is an immediate corollary from Theorem 2.2 in \cite{LZ11}. We provide a proof for the convenience of the reader.

\begin{nnumthm}{1.1}
    Let $\mathfrak A$ be a von Neumann algebra.
    \begin{enumerate}
        \item $\mathfrak A$ is strongly clean if and only if there exists a finite number of mutually orthogonal central projections $P_i$ with sum $I$ such that $P_i \mathfrak A$ is of type $I_{n_i}$, $n_i < \infty$.

        \item $\mathfrak A$ is strongly $*$-clean if and only if $\mathfrak A$ is abelian.
    \end{enumerate}
\end{nnumthm}

\begin{proof}
    By Theorem 6.5.2 in \cite{KR97}, Lemmas \ref{lem:infinite_not_strong_clean}, \ref{lem:shift_tensor_not_strong_clean}, \ref{lem:type_II_1_not_strong_clean}, and \ref{lem:I_strongly_clean}, we have $\FA$ is strongly clean if and only if  there exists a finite number of mutually orthogonal central projections $P_i$ with sum $I$ such that $P_i \mathfrak A$ is of type $I_{n_i}$, $n_i < \infty$. 

We assume that $\FA$ is a non-abelian strongly clean von Neumann algebra. Since $\FA$ is non-commutative, there exist two non-zero mutually orthogonal projections $E_1$ and $E_2$ such that $E_1 \sim E_2$ in $\FA$. Let $T = E_1 + V$, where $V$ is a partial isometry such that $V^*V = E_2$ and $VV^* = E_1$. If $P$ is an projection such that $PT = TP$, then there exists a subprojection $P_1$ of $E_1$ and a subprojection $P_2$ of $I-E_1-E_2$ such that $P = P_1 + V^*P_1V+P_2$. It is clear that $T-P$ is not invertible. Thus $\FA$ is not strongly $*$-clean.
\end{proof}

\begin{remark}
A unital ring $\FR$ is called directly finite if $ab=1$ in $\FR$ implies $ba=1$. In \cite[Question 2]{N99}, Nicholson proposed a question: Is every strongly clean ring directly finite? Let $\FA$ be a finite von Neumann algebra and
$T,S\in\FA$ such that $TS=I$. Then we have
$ker(S)=\{0\}$ and $S$ is lower bounded. Thus $S$ is invertible and $ST=I$. It follows that every finite von Neumann algebra is directly finite. Combining this result with Theorem \ref{thm:char_strongly_clean}.(1), we know that strongly clean von Neumann algebras are directly finite. Therefore, for von Neumann algebras, the answer to Nicholson's question is yes.
\end{remark}

\section{Auxiliary results}\label{sec:aux_results}
This section contains some auxiliary results to be used in later proofs. Some of them are known to experts, we sketch the proofs for the sake of completeness. For every $T \in \CB(\CH)$, we use $\CK(T)$ and $\CR(T)$ to denote the projections onto $\ker(T):=\left \{\xi \in \CH: T\xi = 0 \right\}$ and the closure of the range $\Ran(T):= \left \{T\xi: \xi \in \CH \right\}$ of $T$ (note that $\Ran(T)$ is not necessarily closed), respectively.

Recall that an operator $T$ in a von Neumann algebra $\FA$ is called \textit{finite relative to $\FA$} if there exists a finite projection $E \in \FA$ such that $T \in E \FA E$. And $T \in \FA$ is called a \textit{compact operator relative to $\FA$} (abbreviated as compact operator in $\FA$) if there exists a sequence of operators $T_n$ which are finite relative to $\FA$ such that $\lim_{n\to\infty}\|T- T_n\| = 0$. It is well known that the set of compact operators in $\FA$ is the smallest closed two sided ideal of $\FA$ containing the finite projections of $\FA$ (see \cite{MB68}). In the following proposition, we characterize the operators which can be written as the sum of a scalar and a compact operator in an infinite factor.

\begin{proposition}\label{prop:not_scalar_compact_des}
    Let $\FA$ be an infinite factor and $T \in \FA$. The following statements are equivalent:
    \begin{enumerate}
        \item[(1)] There exists $c \in \mathbb{C}$ such that $T - cI$ is compact relative to $\FA$.
        \item[(2)] $T - W^*TW$ is compact relative to $\FA$ for every unitary $W \in \FA$.
        \item[(3)] $(I-E)TE$ is compact relative to $\FA$ for every infinite projection $E$ such that $E \sim I-E$.
    \end{enumerate}
\end{proposition}

\begin{proof}
    It is clear that (1) implies (2) and (3). By considering the real and imaginary parts of $T$, we only need to show that both (2) and (3) imply (1) under the assumption that $T$ is self-adjoint. Let $H$ be a self-adjoint operator in $\FA$.

    (2) $\Rightarrow$ (1): Assume that $H - W^*HW$ is compact relative to $\FA$ for every unitary $W \in \FA$. We claim that for every $c_1 \in \left [ 0, \|H\| \right]$ and $c_2 < c_1$, the spectral projection $P_1$ of $H$ associated with $\left[c_1, \|H\| \right]$ and the spectral projection $P_2$ of $H$ associated with $\left[-\|H\|, c_2 \right]$ can not both be infinite. Indeed, if $P_1$ and $P_2$ are both infinite, then there exist infinite subprojections $E_1 \leq P_1$ and $E_2 \leq P_2$ such that $E_1 \sim E_2$. Let $W$ be a unitary in $\FA$ such that $W^*E_2W = E_1$. Note that
    \begin{align*}
        \braket{\xi}{(H - W^*HW)\xi} \geq c_1-c_2
    \end{align*}
for every unit vector $\xi \in E_1\CH$. We obtain a contradiction with the fact that $H-W^*HW$ is compact relative to $\FA$. And the claim is proved.

We are now ready to show that there exists $c \in \R$ such that $H - cI$ is compact relative to $\FA$. If $H$ is compact relative to $\FA$, then there is nothing to prove. Assume now that $H$ is not compact relative to $\FA$. Without loss of generality, we may assume that the spectral projection of $H$ associated with $\left [0, \|H\| \right]$ is infinite. 

Let 
$c = \sup \left \{a \in \R: \mbox{the spectral projection of $H$ associated with $[-\|H\|, a]$ is finite} \right \}$. Let $K = H - cI$. If $K = 0$, we are done. If $K \neq 0$, then the spectral projection of $K$ associated with $\left [-\|K\|, -\varepsilon \right ]$ is finite and the spectral projection of $K$ associated with $[-\varepsilon, \varepsilon]$ is infinite for every $0 < \varepsilon < \|K\|$. Since $K - W^*KW$ is compact for every unitary $W \in \FA$, the claim above implies that the spectral projection associated with $\left [\varepsilon, \|K\| \right ]$ is finite. Therefore $K$ is compact relative to $\FA$.

    (3) $\Rightarrow$ (1): Assume that $(I-E)HE$ is compact relative to $\FA$ for every infinite projection $E$ such that $E \sim I-E$. Let $E_0$ be a projection such that $E_0 \sim I-E_0$. We can identify $\FA$ with $E_0\FA E_0 \otimes M_2(\mathbb{C})$ and write $H$ as a matrix of operators
    \begin{align*}
       H = \begin{pmatrix}
           H_{1} & 0\\
           0 & H_{2}
       \end{pmatrix} + K_1
    \end{align*}
where $H_{i} \in E_0\FA E_0$ and $K_1$ is a compact operator in $\FA$. Note that
 \begin{align*}
    \frac{1}{2}\begin{pmatrix}
        I & I \\
        I & -I
    \end{pmatrix}
    \begin{pmatrix}
        H_{1} & 0\\
        0 & W^*H_{2}W
    \end{pmatrix}
    \begin{pmatrix}
        I & I \\
        I & -I
    \end{pmatrix} =
    \frac{1}{2}\begin{pmatrix}
        H_1+W^*H_2W & H_1-W^*H_2W\\
        H_1-W^*H_2W & H_1+W^*H_2W
    \end{pmatrix}
\end{align*}
    for every unitary $W \in E_0\FA E_0$. We have $H_1-W^*H_2W$ is compact relative to $E_0\FA E_0$. In particular $H_1 - H_2$ is compact relative to $\FA$. Therefore, $H_1-W^*H_1W$ is compact relative to $E_0\FA E_0$. Then the equivalence of (1) and (2) implies that there exists $c \in \R$ such that $H - cI$ is compact relative to $\FA$.
\end{proof}

Next, we give some conditions to ensure the finiteness of projections.

\begin{proposition}\label{lem:conjugate_operator}
    Let $T \in \CB(\CH)$ and $a, c \in (0, \infty)$. If $E$ is an idempotent in $\CB(\CH)$ such that $\|TE\| \leq a$ and $c \|\xi\| \leq \|T(I-E)\xi\|$ for every $\xi \in (I-E) \CH$, then the projection $F = I - \CR(T(I-E))$ satisfies the following conditions:
   \begin{enumerate}
       \item $\|T^*F\| \leq a$,
       \item $c \|\beta\| \leq \|T^*(I-F)\beta\|$ for every $\beta \in (I-F)\CH$.
   \end{enumerate}
In particular, if $E$ is the spectral projection of $|T|$ associated with the interval $[0, c]$, then $F$ is the spectral projection of $|T^*|$ associated with the interval $[0,c]$ and $I-E = \CR(T^*(I-F))$.
\end{proposition}

\begin{proof}
  For every $\zeta \in F\CH$,
   \begin{align*}
       \|T^*F \zeta \| = \sup_{\|\xi\|=1} |\braket{\xi}{T^*\zeta}| = \sup_{\|\xi\|=1} |\braket{TE\xi}{\zeta}| \leq \sup_{\|\xi\|=1} \|TE\xi\|\|\zeta\| \leq a\|\zeta\|.
   \end{align*}
We have $\|T^*F\| \leq a$. Since $c \|\xi\| \leq \|T(I-E)\xi\|$ for every $\xi \in (I-E) \CH$, there exists a unique vector $\eta \in (I-E)\CH$ such that $\beta = T\eta$ for every $\beta \in (I-F)\CH$. Note that
    \begin{align*}
        \|\beta\|^2 = \braket{\eta}{T^*\beta} \leq \|\eta\|\|T^*\beta\| \leq \frac{1}{c}\|\beta\|\|T^*\beta\|.
    \end{align*}
We have $c \|\beta\| \leq \|T^*(I-F)\beta\|$. The second part of this proposition can be derived easily from the polar decomposition of bounded operators.
\end{proof}

In the following, Halmos' two projections theorem will be used repeatedly. We refer the reader to \cite{BS10} for a detailed treatment of Halmos' two projections theorem. With Halmos' two projections theorem, we immediately get the following lemma (see also Theorem 1.41 in Chapter V in \cite{T02}).

\begin{lemma}\label{lem:joint_zero_leq}
    Let $E, F$ be two projections in a von Neumann algebra $\FA$. If $F \wedge (I-E) = 0$, then $F \preccurlyeq E$ in $\FA$. In particular $F$ is a finite projection if so is $E$.
\end{lemma}

\begin{corollary}\label{lem:small_finite}
    Let $T$ be an operator in a von Neumann subalgebra $\FA$ of $\CB(\CH)$ and $E$ be the spectral projection of $|T|$ associated with the interval $[0, c]$. Then the following statements hold.
   \begin{enumerate}
       \item[(1)] Let $A \in \FA$ be an operator such that $\|A\| < c$. If $F$ is a projection in $\FA$ such that $\|(T+A)F\| < c - \|A\|$, then $F \preccurlyeq E$.
       \item[(2)] If $c \in (0, 1)$ and $I-E$ is a finite projection, then for every $b \in [0, 1-c)$, the spectral projection $F$ of $|I-T|$ associated with $[0, b]$ is finite.
   \end{enumerate}
\end{corollary}

\begin{proof}
    (1). Let $\xi \in (F \wedge (I-E))\CH$. We claim that $\xi = 0$. Note that $\|T \xi\| \geq c\|\xi\|$ since $E$ is the spectral projection of $|T|$ associated with the interval $[0, c]$. If $\xi \neq 0$, we have
    \begin{align*}
        (c-\|A\|)\|\xi\| \leq \|T\xi\| - \|A \xi\| \leq \|(T+A)\xi\| < (c-\|A\|)\|\xi\|. 
    \end{align*}
     Therefore $\xi = 0$ and $F \wedge (I-E) = 0$. Then \cref{lem:joint_zero_leq} implies that $F \preccurlyeq E$.

    (2). Let $\xi \in (F \wedge E)\CH$. Note that $\|T \xi\| \leq c\|\xi\|$ since $E$ is the spectral projection of $|T|$ associated with the interval $[0, c]$. Therefore,  
    \begin{align*}
        (1-c)\|\xi\| \leq \|\xi\| - \|T\xi\| \leq  \|(I-T)\xi\| \leq b\|\xi\|.
    \end{align*}
    This implies that $\xi = 0$ and $F \wedge E = 0$. By \cref{lem:joint_zero_leq}, we have $F \preccurlyeq I-E$. In particular, $F$ is finite.
\end{proof}

\begin{corollary}\label{lem:range_co_finite}
    Let $T$ be an operator in a von Neumann algebra $\FA$ $(\subseteq \CB(\CH))$. There exists $c > 0$ such that the spectral projection of $|T^*|$ associated with $[0, c]$ is a finite projection if and only if there exists a finite projection $E \in \FA$ such that $(I-E) \CH \subseteq \Ran(T)$.
\end{corollary}

\begin{proof}
    If there exists $c > 0$ such that the spectral projection $E$ of $|T^*|$ associated with $[0, c]$ is finite, then $(I-E) \CH \subseteq \Ran(T)$ by \cref{lem:conjugate_operator}.

    Conversely, assume that $E$ is a finite projection such that $\Ran((I-E)T) = (I-E)\CH$. By the inverse mapping theorem (see Theorem 12.5 in \cite{JC90}), $T$ induces an invertible operator in $\CB((I-\CK((I-E)T))\CH, (I-E)\CH)$. By \cref{lem:conjugate_operator}, there exists $c > 0$ such that $\|T^*\xi\| \geq 2c\|\xi\|$ for every $\xi \in (I-E)\CH$. Let $F$ be the spectral projection of $|T^*|$ associated with $[0, c]$. Note that
    \begin{align*}
        2c \|\xi\| \leq \|T^*\xi\| \leq c\|\xi\|, \quad \forall \xi \in (F \wedge (I-E))\CH.
    \end{align*}
By \cref{lem:joint_zero_leq}, $F \preccurlyeq E$ and $F$ is a finite projection.
\end{proof}

\begin{corollary}\label{cor:compact_range_codim_infinite}
    Let $T \in \FA$. Assume that there exists an infinite projection $E$ in $\FA$ such that $ET$ is compact relative to $\FA$. Then $(I-F)\CH \nsubseteq \Ran(T)$ for every finite projection $F \in \FA$.
\end{corollary}

\begin{proof}
    Since $E$ is an infinite projection and $T^*E$ is compact, the spectral projection of $|T^*|$ associated with $[0, c]$ is infinite for every $c > 0$ by \cref{lem:small_finite}. Then \cref{lem:range_co_finite} implies the result.
\end{proof}

We now provide a criterion for determining when an operator is invertible.

\begin{lemma}\label{lem:two_proj_surjective_cond}
    Let $E, F$ be two projections in $\CB(\CH)$. Then the following are equivalent:
    \begin{enumerate}
        \item $E \wedge F = \{0\}$ and $(E \vee F)\CH = \{\xi + \beta: \xi \in E\CH, \beta \in F\CH\}$.
        \item $E - F$ is an invertible operator in $\CB((E \vee F)\CH)$.
        \item $(E \vee F - F)E$ is an invertible operator in $\CB(E\CH, (E \vee F - F)\CH)$.
    \end{enumerate}
If the above conditions are satisfied, we have
    \begin{align*}
       \left \|((E-F)|_{E\vee F})^{-1} \right \| = \left (1- \|EF\|^2 \right)^{-1/2} = \left\|[(E \vee F-F)E]^{-1} \right \|.
    \end{align*}
\end{lemma}

\begin{proof}
    Without loss of generality, we may assume that $E \vee F = I$. If (2) or (3) holds, we have $E \wedge F = 0$. By Halmos' two projections theorem, we may assume that
   \begin{align*}
      E = \begin{pmatrix}
          I_1 & 0\\
          0 & 0
      \end{pmatrix}
      \oplus \begin{pmatrix}
          I_3 & 0\\
          0 & 0
      \end{pmatrix}, \quad
      F = \begin{pmatrix}
          0 & 0\\
          0 & I_2
      \end{pmatrix}
      \oplus \begin{pmatrix}
          H & \sqrt{H(I_3 - H)}\\
          \sqrt{H(I_3 - H)} & I_3 - H
      \end{pmatrix},
   \end{align*}
where $H$ is a positive contraction such that $\ker(H)= \ker(I_3 - H) = \{0\}$. Then it is clear that (1), (2) or (3) holds if and only if $I_3 - H$ is invertible. Note that $\left \|EF \right \| = \left \|\sqrt{H} \right\|$. We have $\left \|(E-F)^{-1} \right\| = \left (1- \|EF\|^2 \right )^{-1/2} = \left \|[(E \vee F-F)E]^{-1} \right \|$.
\end{proof}

\begin{example}\label{exam:2by2_case}
   Let
   \begin{align*}
      E = \begin{pmatrix}
          H & \sqrt{H(I-H)}\\
          \sqrt{H(I-H)} & I-H
      \end{pmatrix}, \quad
       F = \frac{1}{2}\begin{pmatrix}
           I & -iI\\
           iI & I
       \end{pmatrix}
   \end{align*}
    be two projections in $\CB(\CH) \otimes M_2(\mathbb{C})$, where $H$ is a positive contraction such that $\ker(H) = \ker(I-H) = \{0\}$. It is not hard to check that $E \vee F = I \otimes I_2$, $\|EF\|=\frac{1}{\sqrt{2}}$ and $\left \|(E-F)^{-1} \right \| = \sqrt{2}$.
\end{example}

\begin{theorem}\label{lem:sum_invertible}
    Let $T \in \CB(\CH)$. Then $T$ induces an invertible operator in $\CB(\CH, \CR(T)\CH)$ if and only if there exists a projection $E$ in $\CB(\CH)$ such that the following conditions hold:
   \begin{enumerate}
       \item[(1)] there exist $a_1, a_2 \in (0, \infty)$ such that $\|TE\xi\| \geq a_1 \|E\xi\|$, $\|T(I-E)\xi \| \geq a_2 \|(I-E)\xi\|$ for every $\xi \in \CH$,
       \item[(2)] $\CR(TE) - \CR(T(I-E))$ is invertible in $\CB(\CR(T)\CH)$.
   \end{enumerate}
If the conditions (1) and (2) are satisfied, we have
\begin{align*}
    \frac{1}{\|T\|} \leq \sqrt{\left (1- \|\CR(TE) \CR(T(I-E))\| \right )} \left \|T^{-1} \right \| \leq \frac{1}{\min(a_1, a_2)}.
\end{align*}
\end{theorem}

\begin{proof}
By the inverse mapping theorem and \cref{lem:two_proj_surjective_cond}, the operator $T$ induces an invertible operator in $\CB(\CH, \CR(T) \CH)$ if and only if conditions (1), (2) are satisfied. Let $\xi$ be a unit vector in $\CH$. Note that
\begin{align*}
    \|T\xi\|^2 &\geq \|TE\xi\|^2 + \|T(I-E)\xi\|^2 - 2\|\CR(T(I-E))\CR(TE)\|\|TE\xi\|\|T(I-E)\xi\| \\
    & \geq \inf_{t} \left (\|TE\xi\|^2 + \|T(I-E)\xi\|^2 \right) \left (1-\|\CR(T(I-E))\CR(TE)\|\sin(2t) \right)\\
    & \geq \min \left (a_1^2, a_2^2 \right) \left (1-\|\CR(T(I-E))\CR(TE)\| \right).
\end{align*}
Therefore $\sqrt{(1- \|\CR(TE)\CR(T(I-E))\|)} \left \|T^{-1} \right \| \leq \frac{1}{\min(a_1, a_2)}$.

Without loss of generality, we may assume that
    \begin{align*}
        \CR(TE) = \begin{pmatrix}
       I_1 & 0\\
       0 & 0
   \end{pmatrix} \oplus
    \begin{pmatrix}
        I_3 & 0\\
        0 & 0
    \end{pmatrix} \oplus 0,
    \end{align*}
    \begin{align*}
        \CR(T(I-E)) = \begin{pmatrix}
       0 & 0\\
       0 & I_2
   \end{pmatrix} \oplus \begin{pmatrix}
       H & \sqrt{H(I_3-H)}\\
       \sqrt{H(I_3-H)} & I_3-H
   \end{pmatrix} \oplus 0,
    \end{align*}
where $H$ is a positive contraction such that $\|H\| < 1$ and $\ker(H) = \{0\}$. For every unit vector $\xi:=(0, 0, \xi_1, \xi_2,0)^t \in \CR(T(I-E))\CH$, there exist vectors $\beta \in E \CH$ and $\eta \in (I-E)\CH$ such that
\begin{align*}
    T \beta = \left (0, 0, \frac{\xi_1}{\|\xi_1\|}, 0, 0 \right)^t, \quad T\eta = \xi.
\end{align*}
Note that
\begin{align*}
   \frac{1}{\left \|T^{-1} \right\|^2} \leq \left \|T \left (\frac{\cos t \beta}{\|\beta\|} + \frac{\sin t \eta}{\|\eta\|}\right) \right \|^2 &= \frac{\cos^2 t}{\|\beta\|^2} + \frac{\sin^2 t}{\|\eta\|^2} + \frac{2\cos t \sin t \|\xi_1\|}{\|\beta\|\|\eta\|}\\
    &= \left (\frac{\cos t}{\|\beta\|} + \frac{\sin t}{\|\eta\|} \right)^2 + \frac{2\cos t \sin t \left( \|\xi_1\|-1 \right)}{\|\beta\|\|\eta\|}, \quad \forall t \in \R.
\end{align*}
Choose $t$ such that $\cos t = \frac{\|\beta\|}{\sqrt{\|\beta\|^2 + \|\eta\|^2}}$ and $\sin t = \frac{-\|\eta\|}{\sqrt{\|\beta\|^2 + \|\eta\|^2}}$. Since $\|\beta\| \geq \frac{1}{\|T\|}$ and $\|\eta\| \geq \frac{1}{\|T\|}$, we have
\begin{align*}
    \left (1 - \|\xi_1\| \right) \geq \frac{1}{\left \|T \right\|^2 \left \|T^{-1} \right \|^2}.
\end{align*}
This implies that $\sqrt{1- \|\CR(TE) \CR(T(I-E))\|} \left \|T^{-1} \right\| \geq \frac{1}{\|T\|}$.
\end{proof}

\begin{proposition}\label{lem:idempotent_ran}
    Let $T \in \FA \subseteq \CB(\CH)$. Assume that $E$ is a projection in $\FA$ such that $ETE-E$ is an invertible operator in $E\FA E$. For every $A \in (I-E)\FA E$, there exists an idempotent $P \in \FA$ such that
   \begin{enumerate}
       \item $\CK(P) = I-E$,
       \item $\|P\| \leq 1+ \|TE\|+ \|A\| \left (\|TE\| + 1 \right)$,
       \item $\Ran((T-P)E) = \left \{\xi + A\xi: \xi \in E\CH \right\}$.
   \end{enumerate}
\end{proposition}

\begin{proof}
    Let $P = E+ (I-E)TE - A(ETE-E)$. Since $ETE-E$ is invertible in $E\FA E$, we have $\CK(P)= I-E$, $\|P\| \leq 1+ \|TE\|+ \|A\|(\|TE\| + 1)$, and
    \begin{align*}
        \Ran \left ((T - P)E \right) &= \left \{ (ETE-E)\xi + A(ETE-E)\xi: \xi \in E\CH \right \} = \left \{\xi + A\xi: \xi \in E\CH \right \}.
    \end{align*}
    This completes the proof.
\end{proof}

\section{Finite von Neumann algebras are clean}\label{sec:finite_is_clean}

In this section, we prove that all finite von Neumann algebras are clean.

\begin{lemma}\label{lem:case_easy}
Let $T$ be an operator in a von Neumann algebra $\FA$. Assume that there exist $c > 0$ and a projection $E$ such that
\begin{enumerate}
    \item $ETE -E$ is invertible in $E\FA E$ and $\|T\xi\| \geq c\|\xi\|$ for every $\xi \in (I-E)\CH$,
    \item $E\wedge \CR(T(I-E))\sim (I-E) \wedge (I-\CR(T(I-E)))$ in $\FA$.
\end{enumerate}
Then there exists an idempotent $P$ in $\FA$ such that $\|P\| \leq 2 + 2\|TE\|$, $T-P$ is invertible, and
\begin{align*}
    \left \|(T-P)^{-1} \right \| \leq \frac{2}{\min \left (\left \|(ETE -E)^{-1} \right \|^{-1}, c \right)}.
    \end{align*}
\end{lemma}

\begin{proof}
    Let $F:=I-\CR(T(I-E))$. If $E = 0$, then $F = 0$ and $T$ is invertible and $\left \|T^{-1} \right\| \leq \frac{1}{c}$. Assume that $E \neq 0$ and
    \begin{align*}
        \small{
        E = I_1 \oplus \begin{pmatrix}
            I_2 & 0\\
            0 & 0
        \end{pmatrix}
        \oplus \begin{pmatrix}
           I_3 & 0\\
           0 & 0
        \end{pmatrix} \oplus 0, \quad
        F = I_1 \oplus \begin{pmatrix}
            0 & 0\\
            0 & I_2
        \end{pmatrix} \oplus \begin{pmatrix}
           H & \sqrt{H(I_3-H)}\\
           \sqrt{H(I_3-H)} & I_3-H
       \end{pmatrix} \oplus 0,
   }
    \end{align*}
where $H$ is a positive contraction such that $\ker(H) = \ker(I_3-H) = \{0\}$. By \cref{lem:idempotent_ran}, there exists an idempotent $P$ such that $\CK(P) = I-E$ and
    \begin{align*}
        \CR((T-P)E) = I_1 \oplus
        \begin{pmatrix}
            \frac{1}{2}I_2 & \frac{1}{2}I_2\\
            \frac{1}{2}I_2 & \frac{1}{2}I_2
        \end{pmatrix}
        \oplus \begin{pmatrix}
           \frac{1}{2}I_3 & \frac{-i}{2}I_3\\
           \frac{i}{2}I_3 & \frac{1}{2}I_3
        \end{pmatrix} \oplus 0.
    \end{align*}
By \cref{exam:2by2_case}, it is easy to check that $\left \|\CR(T(I-E)) \CR((T-P)E) \CR(T(I-E)) \right \| = \frac{1}{2}$ and $\CR(T(I-E)) \vee \CR((T-P)E) = I$. By \cref{lem:two_proj_surjective_cond}, $\CR(T(I-E)) - \CR((T-P)E)$ is invertible. Note that $\|(T-P)\xi\| \geq \frac{\|\xi\|}{\|(ETE -E)^{-1}\|}$ for every $\xi \in E\CH$. By \cref{lem:sum_invertible}, we have
    \begin{align*}
        \left \|(T-P)^{-1} \right \| \leq \frac{1}{\sqrt{\left (1- \frac{1}{\sqrt{2}} \right) \min \left (\|(ETE -E)^{-1}\|^{-2}, c^2 \right)}} \leq \frac{2}{\min \left ( \left \|(ETE -E)^{-1} \right \|^{-1}, c \right)}.
    \end{align*}
This completes the proof.
\end{proof}

\begin{example} 
    Let $V \in \CB(\CH)$ be an isometry such that $V$ is not unitary. By Wold decomposition (see, for example, Theorem 1.1 in Chapter I in \cite{NFBK10}), we may assume that $\CH = \CH_0 \oplus (\CH_1 \otimes l^2(\N))$ and $V = V_0 \oplus (I \otimes S)$, where $V_0$ is a unitary on $\CH_0$ and $S$ is the unilateral shift on $l^2(\N)$. Let $E \in \CB(l^2(\N))$ be the projection onto the subspace spanned by $\{e_{2i}\}_{i=1}^{\infty}$, where $\{e_n\}_{n=1}^{\infty}$ is the canonical orthogonal basis of $l^2(\N)$. Recall that $S e_{i} = e_{i+1}$. Then $ESE = 0$, $E \wedge \CR(S(I-E)) = E$, and $(I-E) \wedge (I - \CR(S(I-E))) = I-E$ since $S$ induces a unitary from $(I-E)l^2(\N)$ to $El^2(\N)$. By \cref{lem:case_easy}, $S$ is clean in $\CB(l^2(\N))$. More explicitly, let $P := E + (I-E)(S + S^*)E$, we have $S - P$ is invertible in $\CB(l^2(\N))$. Therefore, $V$ is clean in $\CB(\CH)$.
\end{example}

\cref{thm:finite_is_clean} is an easy corollary of \cref{lem:case_easy}.

\begin{nnumthm}{1.2}
    Let $T$ be an operator in a finite von Neumann algebra $\mathfrak A$. Then there exists an idempotent $P \in \mathfrak A$ such that $T-P$ is invertible and $ \left \|(T-P)^{-1} \right \| \leq 4$. In particular, $\mathfrak A$ is clean.
\end{nnumthm}

\begin{proof}
    If $\|T\| \leq 1/2$, then $T-I$ is invertible and $\left \|(T-I)^{-1} \right \| \leq 2$. Assume that $\|T\| > 1/2$. Let $E$ be the spectral projection of $|T|$ associated with $\left [0, 1/2 \right]$. Since $I-E \sim \CR(T(I-E))$, we have $E \sim I-\CR(T(I-E))$ (see Exercise 6.9.6 in \cite{KR97}). By Halmos' two projections theorem and Exercise 6.9.8 in \cite{KR97}, $E \wedge \CR(T(I-E)) \sim (I-E) \wedge (I- \CR(T(I-E)))$. Note that $\left \|(ETE-E)^{-1} \right \| \leq 2$. By \cref{lem:case_easy}, there exists an idempotent $P \in \FA$ such that $T-P$ is invertible and $\left \|(T-P)^{-1} \right\| \leq 4$.
\end{proof}

\begin{remark}\label{rem:char_almost}
Recall that a $*$-ring is almost $*$-clean if its every ring element can be written as the sum of a projection and a non-zero-divisor (neither a left zero divisor nor a right zero divisor) \cite{LV10}. It is proved in \cite[Corollary 14]{LV10} that finite type I von Neumann algebras are almost $*$-clean. We claim that a von Neumann algebra is almost $*$-clean if and only if it is finite.

Let $\FA$ be a finite von Neumann algebra. An operator in $\FA$ is a non-zero-divisor if and only if its kernal is $\{0\}$. Let $T$ be an operator in $\FA$ such that $\ker(T) \neq \{0\}$. Since $\CK(T) \sim \CK(T^*)$, we may assume that
\begin{align*}
    \small{
    \CK(T) = I_1 \oplus \begin{pmatrix}
        I_2 & 0\\
        0 & 0
    \end{pmatrix}
    \oplus \begin{pmatrix}
        I_3 & 0\\
        0 & 0
    \end{pmatrix} \oplus 0,
}
\end{align*}
\begin{align*}
    \small{
    \CK(T^*) = I_1 \oplus \begin{pmatrix}
        0 & 0\\
        0 & I_2
    \end{pmatrix} \oplus \begin{pmatrix}
        H & \sqrt{H(I_3-H)}\\
        \sqrt{H(I_3-H)} & I_3-H
    \end{pmatrix} \oplus 0,
}
\end{align*}
where $H$ is a positive contraction such that $\ker(H) = \ker(I_3-H) = \{0\}$. Let
\begin{align*}
    P:= I_1 \oplus
        \begin{pmatrix}
            \frac{1}{2}I_2 & \frac{1}{2}I_2\\
            \frac{1}{2}I_2 & \frac{1}{2}I_2
        \end{pmatrix}
        \oplus \begin{pmatrix}
           \frac{1}{2}I_3 & \frac{-i}{2}I_3\\
           \frac{i}{2}I_3 & \frac{1}{2}I_3
        \end{pmatrix} \oplus 0.
\end{align*}
    Note that $\CR(T) = I - \CK(T^*)$. It is not hard to check that $P \wedge \CR(T) = 0$ and $\CK(T) \wedge (I-P) = 0$. Therefore, $\ker(T-P) = \{0\}$. In particular, every finite von Neumann algebra is almost $*$-clean. By the proof of \cref{lem:infinite_not_strong_clean}, we know that properly infinite von Neumann algebras are not almost $*$-clean. In summary, a von Neumann algebra is almost $*$-clean if and only if it is finite.
\end{remark}

\begin{corollary}\label{lem:scalar_plus_compact}
Let $T$ be a compact operator in a von Neumann algebra $\FA$. For every $z \in \mathbb{C}$, there exists an idempotent $P \in \FA$ such that $zI+T - P$ is invertible and $\left \|(zI + T-P)^{-1} \right \| \leq 8$.
\end{corollary}

\begin{proof}
    Since $T$ is compact relative to $\FA$, there exists a finite projection $E \in \FA$ and an operator $A \in E\FA E$ such that $\|T - A\| \leq \frac{1}{8}$. By \cref{thm:finite_is_clean}, there exists an idempotent $P_0 \in E\FA E$ such that $zE + A - P_0$ is invertible in $E\FA E$ and $\left \|(zE + A - P_0)^{-1} \right \| \leq 4$. Let
    \begin{align*}
        P:= \begin{cases}
            P_0, & \mbox{if } |z| \geq \frac{1}{2}, \\
            P_0 + I- E, & \mbox{if } |z| < \frac{1}{2}.
       \end{cases}
    \end{align*}
    Then $zI + A - P$ is invertible and $\left \|(zI + A - P)^{-1} \right \| \leq 4$. Therefore, we have $zI + T-P = \left (zI+A-P \right ) \left [I + (zI+A-P)^{-1}(T-A) \right]$ is invertible and $\|(zI + T-P)^{-1}\| \leq 8$.
\end{proof}

\section{Separable infinite factors are clean}\label{sec:separable_factor_is_clean}

In this section, we prove that all separable infinite factors are clean. From now on, we use $\FA$ to denote an infinite factor (type I$_\infty$, II$_\infty$ or III factor) acting on a separable Hilbert space $\CH$.

\begin{lemma}\label{lem:H_not_easy_case}
Let $T \in \FA$. Assume that $E$ is a projection satisfying the following two conditions:
    \begin{enumerate}
        \item $ETE - E$ is invertible in $E\FA E$ and there exists $c > 0$ such that$\|T\xi\| \geq c \|\xi\|$ for every $\xi \in (I-E)\CH$.
        \item Let $F:=I- \CR(T(I-E))$. There exists $d \in (0, 1)$ such that the spectral projection of $EFE - E\wedge F$ associated with the interval $[0, d]$ is an infinite projection.
    \end{enumerate}
Then there exists an idempotent $P$ such that $T-P$ is invertible and
    \begin{align*}
        \left \|(T-P)^{-1} \right\| \leq \left [ \left (1-\sqrt{\frac{8}{9-d}} \right) \min \left (\left \|(ETE-E)^{-1} \right \|^{-2}, c^2 \right) \right]^{-1/2}.
    \end{align*}
\end{lemma}

\begin{proof}
This result is proved separately for the following two cases.

{\bf Case 1.} $F \wedge (I-E) \preccurlyeq E \wedge (I-F)$: We may assume that
\begin{align*}
    E &= I_1 \oplus
    \begin{pmatrix}
        I_2 & 0\\
        0 & 0
    \end{pmatrix} \oplus
    \begin{pmatrix}
        I_3 & 0 & 0\\
        0 & I_4 & 0\\
        0 & 0 & 0
    \end{pmatrix} \oplus
    \begin{pmatrix}
        I_5 & 0 \\
        0 & 0
    \end{pmatrix} \oplus 0,\\
    F &= I_1 \oplus \begin{pmatrix}
        0 & 0\\
        0 & I_2
    \end{pmatrix} \oplus
    \begin{pmatrix}
        0 & 0 & 0\\
        0 & H_1 & \sqrt{H_1(I_4-H_1)}\\
        0 & \sqrt{H_1(I_4-H_1)} & I_4 -H_1
    \end{pmatrix} \oplus
    \begin{pmatrix}
        H_2 & \sqrt{H_2(I_5-H_2)} \\
        \sqrt{H_2(I_5-H_2)} & I_5 -H_2
    \end{pmatrix} \oplus 0,
\end{align*}
such that $\|H_1\| \leq d$, $0 \oplus
    \begin{psmallmatrix}
        0 & 0\\
        0 & 0
    \end{psmallmatrix} \oplus
    \begin{psmallmatrix}
        0 & 0 & 0\\
        0 & I_4 & 0\\
        0 & 0 & 0
    \end{psmallmatrix} \oplus
    \begin{psmallmatrix}
        0 & 0 \\
        0 & 0
    \end{psmallmatrix} \oplus 0$ is an infinite projection, and $\ker(H_1) = \ker(H_2) = \ker(I_4-H_1) = \ker(I_5-H_2) = \{0\}$. By \cref{lem:idempotent_ran}, there exists an idempotent $P$ such that $\CK(P) = I-E$ and
    \begin{align*}
        \CR((T-P)E) = \CR \left ( I_1 \oplus
\begin{pmatrix}
        I_2 & 0\\
        I_2 & 0
    \end{pmatrix} \oplus
    \begin{pmatrix}
        I_3 & 0 & 0\\
        0 & I_4 & 0\\
        \frac{2}{\sqrt{1-d}}V & \frac{2}{\sqrt{1-d}} W & 0
    \end{pmatrix} \oplus
    \begin{pmatrix}
        I_5 & 0 \\
        iI_5 & 0
    \end{pmatrix} \oplus 0 \right ),
    \end{align*}
where $V$ and $W$ are partial isometries such that $V^*V = I_3$, $W^*W = I_4$ and $VV^* + WW^* = I_4$. Let $\xi:= \left (0, 0, 0, V^*\zeta, W^*\zeta, \frac{2}{\sqrt{1-d}}\zeta, 0, 0, 0 \right)^{t}$ be a unit vector in $\CR((T-P)E)$. Note that
\begin{align*}
    F\xi = \left (0, 0, 0, 0, \sqrt{H_1}\beta, \sqrt{I_4-H_1}\beta, 0, 0, 0 \right)^t,
\end{align*}
where $\beta = \sqrt{H_1}W^*\zeta + \frac{2}{\sqrt{1-d}}\sqrt{I_4-H_1}\zeta$. Since $\|H_1\| \leq d$ and $\|\zeta\| = \sqrt{\frac{1-d}{5-d}}$, we have
\begin{align*}
   \left \|\sqrt{H_1}W^*\zeta + \frac{2}{\sqrt{1-d}}\sqrt{I_4-H_1}\zeta \right \| \geq \|\zeta\|,
\end{align*}
and $\|F\xi\| \geq \sqrt{\frac{1-d}{5-d}}$. This implies that $\Ran(F\CR((T-P)E)) = F\CH$ and $\CR((T-P)E)) \vee (I-F) = I$ (see also \cref{exam:2by2_case}). Note that $\left \|(I-F)\CR((T-P)E) \right \| \leq \sqrt{\frac{4}{5-d}}$. By \cref{lem:two_proj_surjective_cond} and \cref{lem:sum_invertible}, we have
\begin{align*}
    \left \|(T -P)^{-1} \right \| \leq \left [ \left (1-\sqrt{\frac{4}{5-d}}\right )\min \left ( \left \|(ETE-E)^{-1} \right \|^{-2}, c^2 \right ) \right ]^{-1/2}.
\end{align*}

{\bf Case 2.} $E \wedge (I-F) \preccurlyeq F \wedge (I-E)$: We may assume that
\begin{align*}
    E &= I_1 \oplus
    \begin{pmatrix}
        I_2 & 0\\
        0 & 0
    \end{pmatrix} \oplus
    \begin{pmatrix}
        0 & 0 & 0\\
        0 & I_4 & 0\\
        0 & 0 & 0
    \end{pmatrix} \oplus
    \begin{pmatrix}
        I_5 & 0 \\
        0 & 0
    \end{pmatrix} \oplus 0,\\
    F &= I_1 \oplus \begin{pmatrix}
        0 & 0\\
        0 & I_2
    \end{pmatrix} \oplus
    \begin{pmatrix}
        I_3 & 0 & 0\\
        0 & H_1 & \sqrt{H_1(I_4-H_1)}\\
        0 & \sqrt{H_1(I_4-H_1)} & I_4 -H_1
    \end{pmatrix} \oplus
    \begin{pmatrix}
        H_2 & \sqrt{H_2(I_5-H_2)} \\
        \sqrt{H_2(I_5-H_2)} & I_5 -H_2
    \end{pmatrix} \oplus 0,
\end{align*}
such that $\|H_1\| \leq d$, $0 \oplus
    \begin{psmallmatrix}
        0 & 0\\
        0 & 0
    \end{psmallmatrix} \oplus
    \begin{psmallmatrix}
        0 & 0 & 0\\
        0 & I_4 & 0\\
        0 & 0 & 0
    \end{psmallmatrix} \oplus
    \begin{psmallmatrix}
        0 & 0 \\
        0 & 0
    \end{psmallmatrix} \oplus 0$ is an infinite projection, and $\ker(H_1) = \ker(H_2) = \ker(I_4-H_1) = \ker(I_5-H_2) = \{0\}$. By \cref{lem:idempotent_ran}, there exists an idempotent $P$ such that $\CK(P) = I-E$ and
\begin{align*}
    \CR((T-P)E) = \CR \left (I_1 \oplus
\begin{pmatrix}
        I_2 & 0\\
        I_2 & 0
    \end{pmatrix} \oplus
    \begin{pmatrix}
        0 & \frac{2\sqrt{2}}{\sqrt{1-d}}V & 0\\
        0 & I_4 & 0\\
        0 & \frac{2\sqrt{2}}{\sqrt{1-d}} W & 0
    \end{pmatrix} \oplus
    \begin{pmatrix}
        I_5 & 0 \\
        iI_5 & 0
    \end{pmatrix} \oplus 0 \right),
\end{align*}
where $V$ and $W$ are partial isometries such that $VV^* = I_3$, $WW^* = I_4$ and $V^*V + W^*W = I_4$. Let $\xi:= \left (0, 0, 0, \frac{2\sqrt{2}}{\sqrt{1-d}}V\zeta, \zeta, \frac{2\sqrt{2}}{\sqrt{1-d}} W \zeta, 0, 0, 0 \right)^t$ be a unit vector in $\CR((T-P)E)$. We have
\begin{align*}
    F\xi = \left (0,0, 0, \frac{2\sqrt{2}}{\sqrt{1-d}}V\zeta, \sqrt{H_1}\beta, \sqrt{I_4-H_1}\beta, 0, 0, 0 \right)^t,
\end{align*}
where $\beta = \sqrt{H_1}\zeta + \frac{2\sqrt{2}}{\sqrt{1-d}}\sqrt{I_4-H_1}W\zeta$. Since $\|\zeta\|^2 =\|V \zeta\|^2 + \|W\zeta\|^2$, we have $\|\zeta\| = \sqrt{\frac{1-d}{9-d}}$ and
\begin{align*}
    \|F\xi\|^2  \geq
    \begin{cases}
        \left \|\frac{2\sqrt{2}}{\sqrt{1-d}}V\zeta \right \|^2 \geq \|\zeta\|^2, & \text{if }\|V\zeta\|^2 \geq \frac{1}{2}\|\zeta\|^2,\\
        \left \|\sqrt{H_1}\zeta + \frac{2\sqrt{2}}{\sqrt{1-d}}\sqrt{I_4-H_1}W\zeta \right \|^2 \geq \|\zeta\|^2, & \text{if }\|W\zeta\|^2 \geq \frac{1}{2}\|\zeta\|^2. \\
    \end{cases}
\end{align*}
Then it is not hard to see that $(I-F) \vee  \CR((T-P)E) = I$ and $\|(I-F)\CR((T-P)E)\| \leq \sqrt{\frac{8}{9-d}}$. By \cref{lem:two_proj_surjective_cond} and \cref{lem:sum_invertible}, we have
\begin{align*}
    \left \|(T - P)^{-1} \right \| \leq \left [ \left (1-\sqrt{\frac{8}{9-d}} \right ) \min \left ( \left \|(ETE-E)^{-1} \right \|^{-2}, c^2 \right ) \right]^{-1/2}.
\end{align*}
This completes the proof.
\end{proof}

\begin{remark}
Let $T \in \FA$ and $E$ be a projection in $\FA$. The condition (2) in \cref{lem:H_not_easy_case} holds if and only if $E\CR(T(I-E))E - E \wedge \CR(T(I-E))$ is not compact relative to $\FA$. In particular, the condition (2) in \cref{lem:H_not_easy_case} holds if $\FA$ is a type III factor and $E\CR(T(I-E))E - E \wedge \CR(T(I-E)) \neq 0$.
\end{remark}

Combining \cref{lem:case_easy} and \cref{lem:H_not_easy_case}, we obtain the following result.

\begin{lemma}\label{lem:easy_hard}
    Let $T \in \FA$. Assume that there exist $c > 0$ and a projection $E \in \FA$ such that:
   \begin{enumerate}
       \item $ETE - E$ is invertible in $E\FA E$ and $\|T\xi\| \geq c\|\xi\|$ for every $\xi \in (I-E)\CH$,
       \item $E \wedge \CR(T(I-E)) \sim (I-E) \wedge (I- \CR(T(I-E)))$ or $E \CR(T(I-E))E - E \wedge \CR(T(I-E))$ is not compact relative to $\FA$.
   \end{enumerate}
Then every operator $T_1$ similar to $T$ in $\FA$, i.e., $T_1 = V^{-1}TV$ for an invertible operator $V \in \FA$, is clean.
\end{lemma}

\begin{lemma}\label{lem:I_H_compact_12_compact}
    Let $T \in \FA$. If $E$ is a projection in $\FA$ such that $E\CR(T(I-E))E - E \wedge \CR(T(I-E))$ is compact relative to $\FA$, then $[E- E \wedge \CR(T(I-E))]T(I-E)$ is compact relative to $\FA$.
\end{lemma}

\begin{proof}
    By Halmos' two projections theorem,
    \begin{align*}
        \Ran \left (\left [E- E \wedge \CR(T(I-E)) \right]T(I-E) \right) \subset \Ran\left(\sqrt{E\CR(T(I-E))E - E \wedge \CR(T(I-E))}\right).
    \end{align*}
    Since $E\CR(T(I-E))E - E \wedge \CR(T(I-E))$ is compact relative to $\FA$, $\left [E- E \wedge \CR(T(I-E)) \right]T(I-E)$ is also compact relative to $\FA$.
\end{proof}

\begin{lemma}\label{lem:two_condition_not_satisfied}
   Let $T \in \FA$ and $E$ be a projection in $\FA$ satisfying the following conditions:
    \begin{enumerate}
        \item $E \sim I-E$,
        \item there exists $c > 0$ such that $\|T\xi\| \geq c\|\xi\|$ for every $\xi \in (I-E)\CH$,
        \item $E \wedge \CR(T(I-E)) \nsim (I-E) \wedge (I-\CR(T(I-E)))$,
        \item $E\CR(T(I-E))E - E \wedge \CR(T(I-E))$ is compact relative to $\FA$.
    \end{enumerate}
    If $ET(I-E)$ is not compact relative to $\FA$, then $E \wedge \CR(T(I-E))$ is an infinite projection and there exists a finite subprojection $F$ of $I-E$ such that $(I-E-F)\CH \subset \Ran((I-E)T(I-E))$. Furthermore, if there exists a finite subprojection $E_0$ of $E$ such that $(E-E_0)\CH \subseteq \Ran(ET(I-E))$, then there exists $a > 0$ such that the spectral projection of $|T^*|$ associated with $[0, a]$ is finite and the spectral projection of $|T|$ associated with $[0, \varepsilon]$ is infinite for every $\varepsilon > 0$.
\end{lemma}

\begin{proof}
    By \cref{lem:I_H_compact_12_compact}, $[E - E \wedge \CR(T(I-E))]T(I-E)$ is compact relative to $\FA$. Since $ET(I-E)$ is not compact relative to $\FA$, we have $E \wedge \CR(T(I-E))$ is infinite. Thus $(I-E) \wedge (I- \CR(T(I-E)))$ is a finite projection by condition (3) and Corollary 6.3.5 in \cite{KR97}. Let $P$ be the subprojection of $I-E$ such that $\Ran(TP) = [\CR(T(I-E)) - E \wedge \CR(T(I-E))]\CH$. Since $[E - E \wedge \CR(T(I-E))]TP$ is compact relative to $\FA$, there exists a finite subprojection $P_1$ of $P$ such that
\begin{align*}
    \left \|(I-E)T\beta \right \| \geq \frac{c}{2}\|\beta\|, \quad \forall \beta \in (P-P_1)\CH.
\end{align*}
In particular, $\Ran((I-E)T(P-P_1))$ is closed. Let $F = (I-E)-\CR((I-E)T(P-P_1))$. Since $P_1$ and
\begin{align*}
    I-E - \CR((I-E)T(I-E)) &= (I-E) \wedge (I-\CR(T(I-E)))\\
    & =(I-E) \wedge (I-\CR(TP)) = (I-E) - \CR((I-E)TP)
\end{align*}
are finite projections, $F$ is a finite projection.

    Assume that there exists a finite subprojection $E_0$ of $E$ such that $(E-E_0)\CH \subseteq \Ran(ET(I-E))$. By \cref{cor:compact_range_codim_infinite}, we have $E - E \wedge \CR(T(I-E))$ is a finite projection. Then it is not hard to see that $(I-E) - (I-E) \wedge \CR(T(I-E))$ is a finite projection. Thus
    \begin{align*}
        I- \left [E \wedge \CR(T(I-E)) + (I-E) \wedge \CR(T(I-E)) \right]
    \end{align*}
    is a finite projection. Since $[E \wedge \CR(T(I-E)) + (I-E) \wedge \CR(T(I-E))]\CH \subset \Ran(T)$, there exists $a > 0$ such that the spectral projection of $|T^*|$ associated with $[0, a]$ is finite by \cref{lem:range_co_finite}.

    To prove that the spectral projection of $|T|$ associated with $[0, \varepsilon]$ is infinite for every $\varepsilon > 0$, we only need to show that there exists an infinite projection $Q$ such that $\|TQ\| < \varepsilon$. If $TE$ is compact relative to $\FA$, then there exists an infinite projection $Q \leq E$ such that $\|TQ\| < \varepsilon$. Assume that $TE$ is not compact relative to $\FA$. Then there exists an infinite projection $Q_0$ such that $Q_0 \CH \subseteq \Ran(TE)$. Since $I- [E \wedge \CR(T(I-E)) + (I-E) \wedge \CR(T(I-E)]$ is a finite projection, $Q_0 \wedge \CR(T(I-E))$ is an infinite projection. This implies that $\CK(T)$ is an infinite projection and the lemma is proved.
\end{proof}

\begin{lemma}\label{lem:case_I}
    Let $T \in \FA$. Assume that $T - cI$ is not compact relative to $\FA$ for every $c \in \mathbb{C}$ and there exists $a \in (0, 1)$ such that the spectral projection of $|T|$ associated with $[0, a]$ is finite, then there exists an invertible operator $V$ and a projection $E$ in $\FA$ satisfying the following conditions:
   \begin{enumerate}
       \item $E \sim I-E$,
       \item $\left \|EV^{-1}TVE \right \| < 1$,
       \item there exists $c > 0$ such that $\left \|V^{-1}TV\xi \right \| \geq c \left \|\xi \right \|$ for every $\xi \in (I-E) \CH$,
       \item there exists a finite subprojection $F$ of $E$ such that $(E-F)\CH \subset \Ran \left (EV^{-1}TV(I-E) \right)$.
   \end{enumerate}
\end{lemma}

\begin{proof}
    By \cref{prop:not_scalar_compact_des}, there exists a projection $E_0$ such that $E_0 \sim I-E_0$ and $(I-E_0)T^*E_0$ is not compact relative to $\FA$. Since $(I-E_0)T^*E_0$ is not compact relative to $\FA$,  there exists $b > 0$ such that the spectral projection $E_1$ of $|(I-E_0)T^*E_0|$ associated with $\left [b, \|(I-E_0)T^*E_0\| \right ]$ is an infinite subprojection of $E_0$. Note that there exists an operator $A$ in $E_0\FA (I-E_0)$ such that $A(I-E_0)T^*E_0 = E_0T^*E_1$. Let
    \begin{align*}
        V: = \left (I-A^* \right ) \left [E_1 + 2\|(I+A^*)T(I-A^*)\|(I-E_1)\right].
    \end{align*}
    It is not hard to see that $E_1 V^{-1}TV E_1 = 0$, $\left \|(I-E_1)V^{-1}TV E_1 \right \| \leq \frac{1}{2}$, and
    $$\left \|(I-E_1) V^*T^*(V^*)^{-1} E_1 \xi \right \| \geq 2b \left \|(I+A^*)T(I-A^*) \right \|\|\xi\|$$ for every $\xi \in E_1 \CH$. In particular, we have $\Ran \left (E_1 V^{-1}TV (I-E_1) \right ) = E_1 \CH$.

    By \cref{lem:range_co_finite}, there exists $c \in (0, \frac{1}{2})$ such that the spectral projection of $|V^{-1}TV|$ associated with $[0, 2c]$ is finite. Let $E$ be the spectral projection of $|V^{-1}TV(I-E_1)|$ associated with $[0, c]$. Note that $E \geq E_1$. Since $\left \|V^{-1}TV(E - E_1) \right \|\leq c < 2c$, $E - E_1$ is a finite subprojection of $I-E_1$ by \cref{lem:small_finite}. By Corollary 6.3.5 in \cite{KR97}, $E \sim I-E$. Note that
    \begin{align*}
        \left \|E V^{-1}TVE \right\| \leq \left \|(E-E_1)V^{-1}TVE_1 \right \| + \left \|EV^{-1}TV(E-E_1) \right\| < 1,
    \end{align*}
    and $\left \|V^{-1}TV(I-E)\xi \right\| \geq c\|\xi\|$ for every $\xi \in (I-E)\CH$. Since $EV^{-1}TV(I-E) - E_1V^{-1}TV(I-E_1)$ is finite relative to $\FA$, there exists a finite subprojection $F$ of $E$ such that $(E-F) \CH \subset \Ran \left (EV^{-1}TV(I-E) \right)$.
\end{proof}

\begin{lemma}\label{lem:12_not_compact_11_small}
   Let $T \in \FA$. Assume that the spectral projection of $|T|$ associated with $[0, c]$ and its complement are both infinite for every $c \in (0, 1)$. For every $a > 0$, there exist an invertible operator $V$ and a projection $E$ in $\FA$ satisfying the following conditions:
    \begin{enumerate}
        \item $E \sim I-E$,
        \item $\left \|EV^{-1}TVE \right\| \leq \frac{3}{4}$,
        \item For every $\varepsilon > 0$, there exists an infinite subprojection $E_0$ of $E$ such that $E_0 \sim E-E_0$ and $\left \|V^{-1}TV\xi \right \| \leq \varepsilon \|\xi\|$ for every $\xi \in E_0\CH$,
        \item $\left \|V^{-1}TV\xi \right \| \geq \frac{1}{4}\|\xi\|$ for every $\xi \in (I-E)\CH$,
        \item $EV^{-1}TV(I-E)$ is not compact relative to $\FA$,
        \item $\left \|V-I \right \| \leq a$.
    \end{enumerate}
\end{lemma}

\begin{proof}
    Let $E$ be the spectral projection of $|T|$ associated with $\left [0, \frac{1}{2} \right]$. If $ET(I-E)$ is not compact relative to $\FA$, then we are done.

    Assume that $ET(I-E)$ is compact relative to $\FA$. Let $W \in E\FA(I-E)$ be a partial isometry such that $WW^* = E$ and $W^*W = I-E$. Since $\|T\xi\| \geq \frac{1}{2}\|\xi\|$ for every $\xi \in (I-E)\CH$, there exists a finite projection $E_1$ of $I-E$ such that $\|(I-E)T\beta\| \geq \frac{1}{4}\|\beta\|$ for every $\beta \in (I-E-E_1)\CH$. Let $E_2$ be an infinite subprojection of $E$ such that $\|TE_2\| \leq \frac{1}{8}$. Note that $W^*E_2W \wedge (I-E-E_1)$ is an infinite projection and
\begin{align*}
    \left \| \left [WT(I-E) - ETW \right ]\zeta \right \| \geq \frac{1}{8}\|\zeta\|, \quad \forall \zeta \in  \left [W^*E_2W \wedge (I-E-E_1) \right]\CH.
\end{align*}
In particular, $WT(I-E) - ETW$ is not compact relative to $\FA$.

For every $t > 0$, $I-tW$ is invertible. Note that
    \begin{align*}
        E(I - tW)T(I+tW)(I-E) = ET(I-E) - t \left (WT(I-E) - ETW \right) - t^2WTW
    \end{align*}
    is compact relative to $\FA$ if and only if $(WT(I-E) - ETW) + tWTW$ is compact relative to $\FA$. Since $WT(I-E) - ETW$ is not compact relative to $\FA$, there exists $t_0 \in \left (0, \min \left (a, \frac{1}{16\|T\|} \right ) \right)$ such that $E (I-t_0W)T(I+t_0W)(I-E)$ is not compact relative to $\FA$. It is not hard to see that $I+t_0W$ and $E$ satisfy the conditions (1)-(6).
\end{proof}

We are now ready to prove \cref{thm:separable_factor_is_clean}.

\begin{nnumthm}{1.3}
   Every separable infinite factor is clean.
\end{nnumthm}

\begin{proof}
Let $T \in \FA$. In order to show that $T$ is clean, we only need to consider the following two cases:
\begin{enumerate}
    \item[1] There exists $c \in (0, 1)$ such that either the spectral projection $E_c$ of $|T|$ associated with $[0, c]$ is finite or its complement $I-E_c$ is finite.
    \item[2] For every $c \in (0, 1)$, the spectral projections of $|T|$ and $|T^*|$ associated with $[0, c]$ and their complements are infinite.
\end{enumerate}

    {\bf Case 1}: Note that $T$ is clean if and only if $I-T$ is clean. By \cref{lem:small_finite}, we may assume that there exists $c \in (0, 1)$ such that the spectral projection of $|T|$ associated with $[0, c]$ is finite. If there exists $z \in \mathbb{C}$ such that $T-zI$ is compact relative to $\FA$, then $T$ is clean by \cref{lem:scalar_plus_compact}. Assume that $T-zI$ is not compact relative to $\FA$ for every $z \in \mathbb{C}$. There exist an invertible operator $V$ and a projection $E$ in $\FA$ satisfying conditions (1)-(4) in \cref{lem:case_I}. By condition (4) in \cref{lem:case_I} and \cref{lem:two_condition_not_satisfied}, we have either $E \wedge \CR \left (V^{-1}TV(I-E) \right ) \sim (I-E) \wedge \left (I- \CR \left (V^{-1}TV(I-E) \right ) \right )$ or $E\CR \left (V^{-1}TV(I-E) \right)E - E \wedge \CR \left (V^{-1}TV(I-E) \right)$ is not compact. Then $T$ is clean by \cref{lem:easy_hard}.

    {\bf Case 2}: Let $\{E_{ij}\}$ be the canonical matrix units of $M_2(\mathbb{C})$. By \cref{lem:12_not_compact_11_small}, we can assume that $\CH = \CH_1 \otimes \mathbb{C}^2$, $\FA = \FA_1 \otimes M_2(\mathbb{C})$, and
    \begin{align*}
       T = \begin{pmatrix}
           T_{11} & T_{12}\\
           T_{21} & T_{22}
       \end{pmatrix}
    \end{align*}
    such that
    \begin{enumerate}
        \item $\|(I \otimes E_{11})T(I \otimes E_{11})\| \leq \frac{3}{4}$,
        \item for every $\varepsilon > 0$, there exists a projection $F \in \FA_1$ such that $F \sim I- F$ in $\FA_1$ and $\|T_{11}F\| \leq \varepsilon$, $\|T_{21}F\| \leq \varepsilon$,
        \item $\|T\xi\| \geq \frac{1}{4}\|\xi\|$ for every $\xi \in (I \otimes E_{22})\CH$,
        \item $T_{12}$ is not compact relative to $\FA_1$.
    \end{enumerate}
By \cref{lem:easy_hard}, we only need to show that $T$ is clean under the assumption that
\begin{align*}
  (I \otimes E_{11}) \wedge \CR \left (T(I \otimes E_{22}) \right) \nsim I \otimes E_{22} \wedge (I \otimes I_2 - \CR(T(I \otimes E_{22})))
\end{align*}
and $(I \otimes E_{11}) \CR(T(I \otimes E_{22})) (I \otimes E_{11}) - (I \otimes E_{11}) \wedge \CR(T(I \otimes E_{22}))$ is compact relative to $\FA$.

Let $E_1$ and $F_1$ be projections in $\FA_1$ such that
\begin{align*}
E_1 \otimes E_{11} = (I \otimes E_{11}) \wedge \CR(T(I \otimes E_{22})), \quad
    \CR(T(F_1 \otimes E_{22})) = \CR(T(I\otimes E_{22})) - E_1 \otimes E_{11}.
\end{align*}
    By \cref{lem:two_condition_not_satisfied}, $F_1 \sim I-F_1$ and there exists a finite projection $F_2 \in \mathfrak{A}_1$ such that $I-F_2 \subset \Ran(T_{22} F_1)$. By \cref{lem:range_co_finite}, there exists $c > 0$ such that the spectral projection of $|F_1T_{22}^*|$ associated with $[0, c]$ is finite.

By condition (2), there exists a unitary operator $W \in \FA_1$ such that
\begin{align*}
    \|T_{11}WF_1\| \leq \frac{c}{4},\quad  \|T_{21}WF_1\| \leq \frac{c}{4}.
\end{align*}
Let $t \in \left (0, \frac{1}{32\|T\|} \right)$ and $V:=I \otimes I_2 + tW \otimes E_{12}$. Note that
\begin{align*}
    V^{-1}TV  = \begin{pmatrix}
       T_{11} - tW T_{21} & T_{12} - t[WT_{22} - T_{11}W + tWT_{21}W]\\
        T_{21} & T_{22} + tT_{21}W
   \end{pmatrix}.
\end{align*}
By \cref{lem:small_finite}, the spectral projection of $|(WT_{22}F_1 - T_{11}WF_1 + tWT_{21}WF_1)^*|$ associated with $[0, \frac{c}{4}]$ is a finite projection. Note that $T_{12}F_1 = (I-E_1)T_{12}F_1$ is compact relative to $\FA_1$. There exists a finite projection $E_2$ in $\FA_1$ such that $(I-E_2) \CH_1 \subset \Ran(T_{12} - t(WT_{22} - T_{11}W + tWT_{21}W))$. It is not hard to check that $\left \|T_{11} - tWT_{21} \right \| \leq \frac{7}{8}$ and $\left \|V^{-1}TV\xi \right\| \geq \frac{1}{8}\|\xi\|$ for every $\xi \in \left (I \otimes E_{22} \right)\CH$. Then $T$ is clean by \cref{lem:easy_hard} and \cref{lem:two_condition_not_satisfied}.
\end{proof}

The following corollary is an immediate consequence of \cref{thm:separable_factor_is_clean} since the homomorphic image of a clean ring is clean.

\begin{corollary}
    The Calkin algebra on a separable infinite-dimensional Hilbert space $\CH$, i.e., the quotient of $\CB(\CH)$ by the ideal of compact operators, is clean.
\end{corollary}

\section{Questions and remarks}\label{sec:concluding remarks}

We proved in Section \ref{sec:finite_is_clean} that all finite von Neumann algebras are clean. It is natural to consider $*$-cleanness of finite von Neumann algebras. It claimed in \cite[Proposition 15]{LV10} that $\oplus_{n=1}^{\infty} M_n(\mathbb{C})$ is $*$-clean. However, in the proof of \cite[Proposition 15]{LV10}, the author used the fact that the direct sum $\oplus_{n=1}^\infty \FA_n$ of a family of finite dimensional C$^*$-algebras $\{\FA_n\}_{n=1}^{\infty}$ is $*$-isomorphic to the direct product $\prod_{n=1}^{\infty} \FA_n$. The proof is invalid since $\prod_{n=1}^{\infty} \FA_n$ is not even a C$^*$-algebra in general. Therefore, we can not use \cite[Proposition 15]{LV10} to deduce that $\oplus_{n=1}^{\infty} M_n(\mathbb{C})$ is $*$-clean. Based on this fact, we propose the following question.

\begin{question}\label{finite von Neumann algebra *-clean}
    Is the finite von Neumann algebra $\oplus_{n=1}^{\infty} M_n(\mathbb{C})$ $*$-clean? More generally, are all finite von Neumann algebras $*$-clean?
\end{question}

\begin{remark}\label{rem:finite von Neumann algebra of type $I_n$ $*$-clean}
    It has been proved in \cite{LV10} (see also \cite{HN01}) that if a ring $\FR$ is $*$-clean, then $M_n(\FR)$ is also $*$-clean. Since every von Neumann algebra of type $I_n$ has the form $\mathcal{A} \otimes  M_n(\mathbb{C})$, where $\mathcal{A}$ is an abelian von Neumann algebra, we know that finite direct sums of type $I_n$ von Neumann algebras are $*$-clean.

    It is also not hard to show that every operator $T$ with closed range in a finite von Neumann algebra acting on a Hilbert space $\mathcal{H}$ is $*$-clean. Indeed, let $E=I-\CK(T)$ and $F=\CR(T)$. Since $E \sim F$, we may assume that
\begin{align*}
    E=I_1\oplus
\begin{pmatrix}
I_2 &  \\
&0\\
\end{pmatrix}
\oplus
\begin{pmatrix}
I_3 &  \\
&0\\
\end{pmatrix}
\oplus 0, \  \
    F=I_1\oplus
\begin{pmatrix}
0 &  \\
&I_2\\
\end{pmatrix}
\oplus
\begin{pmatrix}
H & \sqrt{H(I_3-H)}  \\
\sqrt{H(I_3-H)}&I_3-H
\end{pmatrix}
\oplus 0.
\end{align*}
Let
\begin{equation*}
P=
0\oplus
\begin{pmatrix}
\frac{1}{2}I_2 & \frac{1}{2}I_2 \\
\frac{1}{2}I_2 & \frac{1}{2}I_2\\
\end{pmatrix}
\oplus
\begin{pmatrix}
\frac{1}{2}I_3 & \frac{i}{2}I_3  \\
-\frac{i}{2}I_3 & \frac{1}{2}I_3
\end{pmatrix}
\oplus I_4.
\end{equation*}
Note that $\Ran(P(I-E)) = \Ran(P)$. By \cref{exam:2by2_case}, it is not hard to see that $\Ran(T-P) = \mathcal{H}$. This implies that $T-P$ is invertible.
\end{remark}

\begin{question}
    Given a von Neumann algebra $\mathfrak A$, let
    \begin{align*}
        \mathcal{S} := \{c\geq 1: \mbox{$\forall T \in \mathfrak A$, there exists an idempotent $P$ such that $T-P$ is invertible and $\|(T-P)^{-1}\| \leq c$}\}.
    \end{align*}
Is the set $\mathcal{S}$ non-empty? By \cref{thm:finite_is_clean}, we know that $4 \in \mathcal{S}$ if $\mathfrak A$ is a finite von Neumann algebra. Then we can ask what is the infimum of the set $\mathcal{S}$, providing it is not empty.
\end{question}

\begin{question}
    Does \cref{thm:separable_factor_is_clean} remain true for arbitrary infinite factors? More generally, are all von Neumann algebras clean?
\end{question}

\begin{question}\label{que:que_3}
    Does there exist a clean C$^*$-algebra $\FA$ such that the C$^*$-algebra $\FA \otimes l^{\infty}(\N)$ is not clean?
\end{question}

\begin{remark}
    Let $\FA_n = \{zI_n + A: \mbox{$z \in \mathbb{C}$, $A$ is a strictly upper triangular matrix in $M_n(\mathbb{C})$}\}$. It is clear that $\FA_n$ is a strongly clean Banach algebra. By the proof of \cref{lem:shift_tensor_not_strong_clean}, the element $\oplus_{n=1}^{\infty} V_n$ in $\oplus_{n=1}^{\infty}\FA_n$ is not clean, where $V_n$ is the upper shift matrix in $M_n(\mathbb{C})$. Thus the Banach algebra $\oplus_{n=1}^{\infty}\FA_n$ is not clean. This simple example seems to indicate that the answer to \cref{que:que_3} should be yes.
\end{remark}

\noindent\textit{Acknowledgement.}
The authors would like to thank the anonymous referee for constructive criticisms and valuable comments.

\bibliography{Bib}
\bibliographystyle{amsplain}

\end{document}